\documentclass[aap,seceqn,MSNbibl,citesort,dvips]{arximspdf}

\doi{10.1214/09-AAP628}
\volume{20}
\issue{2}
\pubyear{2010}
\firstpage{593}
\lastpage{639}

\makeatletter

\DeclareMathAlphabet\mathcaligr{OMS}{cmsy}{m}{n}

\newtheorem{theorem}{Theorem}[section]
\newtheorem{prop}[theorem]{Proposition}
\newtheorem{cor}[theorem]{Corollary}
\newtheorem{lem}[theorem]{Lemma}
\newproclaim{defi}[theorem]{Definition}

\makeatother

\begin{document}
\begin{frontmatter}

\title{Interacting Markov chain Monte Carlo methods for solving nonlinear
measure-valued equations\thanksref{T1}}
\runtitle{Interacting Markov chain Monte Carlo methods}

\thankstext{T1}{Supported in part by the French national programs
``S\'{e}curit\'{e} et
Informatique'' and ``Syst\`{e}mes Complexes et Mod\'{e}lisation
Math\'{e}matique'' under
projects NEBBIANO ANR-06-SETI-009 and \mbox{VIROSCOPY} ANR-08-SYSC-016.}

\begin{aug}
\author[A]{\fnms{Pierre} \snm{Del Moral}\corref{}\ead[label=e1]{Pierre.Del-Moral@inria.fr}} and
\author[B]{\fnms{Arnaud} \snm{Doucet}\ead[label=e2]{Arnaud@stat.ubc.ca}}
\runauthor{P. Del Moral and A. Doucet}
\affiliation{INRIA and Universit\'{e} Bordeaux and University of British
Columbia and The~Institute~of~Statistical~Mathematics}
\address[A]{Centre INRIA Bordeaux et Sud-Ouest\\
\quad and Institut de Math\'{e}matiques de Bordeaux\\
Universit\'{e} Bordeaux\\
351 cours de la Lib\'{e}ration\\
33405 Talence cedex\\
France\\
\printead{e1}} 
\address[B]{Department of Statistics\\
\quad and Department of Computer Science\\
University of British Columbia\\
333-6356 Agricultural Road\\
Vancouver, BC, V6T 1Z2\\
Canada\\
and\\
The Institute of Statistical Mathematics\\
4-6-7 Minami-Azabu\\
Minato-ku, Tokyo 106-8569\\
Japan\\
\printead{e2}}
\end{aug}

\received{\smonth{2} \syear{2008}}
\revised{\smonth{7} \syear{2009}}

%
\begin{abstract}
We present a new class of interacting Markov chain Monte Carlo algorithms
for solving numerically discrete-time measure-valued equations. The
associated stochastic processes belong to the class of self-interacting
Markov chains. In contrast to traditional Markov chains, their time
evolutions depend on the occupation measure of their past values. This
general methodology allows us to provide a natural way to sample from a
sequence of target probability measures of increasing complexity. We develop
an original theoretical analysis to analyze the behavior of these iterative
algorithms which relies on measure-valued processes and semigroup
techniques. We establish a variety of convergence results including
exponential estimates and a uniform convergence theorem with respect to the
number of target distributions. We also illustrate these algorithms in the
context of Feynman--Kac distribution flows.
\end{abstract}

%

\end{frontmatter}

\section{Introduction}\label{sec1}

\subsection{Nonlinear measure-valued processes}

Let $(S^{(l)},\mathcaligr{S}^{(l)})_{l\geq0}$ be a sequence of measurable
spaces. For every $l\geq0$ we denote by $\mathcaligr{P}(S^{(l)})$ the
set of
probability measures on $S^{(l)}$. Suppose we have a sequence of probability
measures $\pi^{(l)}\in\mathcaligr{P}(S^{(l)})$ where $\pi^{(0)}$ is known
and we have for $l\geq1$ the following nonlinear measure-valued equations
%
\begin{equation}\label{phi}
\pi^{(l)}=\Phi_{l}\bigl(\pi^{(l-1)}\bigr)
\end{equation}
for some mappings $\Phi_{l}\dvtx\mathcaligr{P}(S^{(l-1)})\rightarrow
\mathcaligr{P}%
(S^{(l)})$. Except in some particular situations, these measure-valued
equations do not admit an analytic solution.

Being able to solve these equations numerically has numerous
applications in
nonlinear filtering, global optimization, Bayesian statistics and
physics as
it would allow us to approximate any sequence of fixed ``target'' probability
distributions $(\pi^{(l)})_{l\geq0}$. For example, in a nonlinear filtering
framework $\pi^{(l)}$ corresponds to the posterior distribution of the state
of an unobserved dynamic model at time $l$ given the observations collected
from time $0$ to time $l$. In an optimization framework, $\pi^{(l)}$ could
correspond to a sequence of annealed versions of a distribution $\pi$ that
we are interested in maximizing. In both cases, $\Phi_{l}$ is a Feynman--Kac
transformation \cite{fk}.

In recent years, there has been considerable interest in the
development of
interacting particle interpretations of measure-valued equations of the form
(\ref{phi}) which we briefly review here.

\subsection{Interacting particle methods}
\label{ips}

The central idea of interacting particle methods is to construct a Markov
chain $X^{(l)}=(X_{p}^{(l)})_{1\leq p\leq N}$ taking values in the product
spaces $(S^{(l)})^{N}$ so that the empirical measure $\pi
_{N}^{(l)}:=\frac{1%
}{N}\sum_{p=1}^{N}\delta_{X_{p}^{(l)}}$ approximates $\pi^{(l)}$ as $
N\uparrow\infty$. In the simpler version, we construct inductively $%
X^{(l)}=(X_{p}^{(l)})_{1\leq p\leq N}$ by sampling $N$ independent random
variables with common law $\Phi_{l}(\pi_{N}^{(l-1)})$. The rationale
behind this is that the resulting particle measure $\pi_{N}^{(l)}$ should
be a good approximation of $\pi^{(l)}$ as long as $\pi_{N}^{(l-1)}$
is a
good approximation of $\pi^{(l-1)}$. More formally, $X^{(l)}$ is an $%
(S^{(l)})^{N}$-valued Markov chain with elementary transitions given by the
following formula:
%
\begin{equation}\quad
\mathbb{P} \bigl( \bigl(X_{1}^{(l)},\ldots,X_{N}^{(l)}\bigr)\in
dx | X^{(l-1)} \bigr)
=\prod_{p=1}^{N}\Phi_{l} \biggl( \frac{1}{N} \sum_{1\leq q\leq
N} \delta
_{X_{q}^{(l-1)}} \biggr) (dx_{p}),
\end{equation}
where $dx=d(x_{1},\ldots,x_{N})=dx_{1}\times\cdots\times dx_{N}$ stands
for an infinitesimal neighborhood of a point in the product space $%
(S^{(l)})^{N}$.

For Feynman--Kac transformations, these interacting particle models
have been
extensively studied and they are sometimes referred to as sequential Monte
Carlo methods, particle filters and population Monte Carlo methods;
see \cite{fk,arnaud} for a review of the literature. In this context, the convergence
analysis of these particle algorithms is now well understood. A variety of
theoretical results are available, including sharp propagations of chaos
properties, fluctuations and large deviations theorems, as well as uniform
convergence results with respect to the level index $l$.

These interacting particle methods suffer from two serious limitations.
First, when the mapping $\Phi_{l}$ is complex, it may be impossible to
generate independent draws from it. Second, it is typically impossible to
determine beforehand the number of particles necessary to achieve a fixed
precision for a given application and users usually have to perform multiple
runs for an increasing number of particles until stabilization of the Monte
Carlo estimates is observed. Markov chain Monte Carlo (MCMC) methods appear
as a natural way to solve these two problems \cite{robert2004}. However,
standard MCMC methods do not apply in this context as we have a
sequence of
target distributions defined on different spaces and the normalizing
constants of these distributions are typically unknown.

\subsection{Self-interacting Markov chains}

We propose here a new class of interacting MCMC methods (i-MCMC) to solve
these nonlinear measure-valued equations numerically. These i-MCMC methods
can be described as adaptive and dynamic simulation algorithms which take
advantage of the information carried by the past history to increase the
quality of the next sequence of samples. Moreover, in contrast to
interacting particle methods, these stochastic algorithms can increase the
precision and performance of the numerical approximations iteratively.

The origins of i-MCMC methods can be traced back to a pair of
articles \cite{dm1,dm2} presented by the first author in collaboration with Laurent Miclo.
These studies are concerned with biology-inspired self-interacting Markov
chain (SIMC) models with applications to genetic type algorithms
involving a
competition between a reinforcement mechanism and a potential
function \cite{dm1,dm2}. These ideas have been extended to the MCMC methodology in the
joint articles of the authors with Christophe Andrieu and Ajay
Jasra \cite{ajay}, as well as in the more recent article of the authors with Anthony
Brockwell \cite{brockwell}. Related ideas have also appeared in
computational chemistry \cite{lyman2006} and statistics \cite{kou2006}.

In the present article, we design a new general class of i-MCMC methods.
Roughly speaking, these algorithms proceed as follows. At level $l=0$
we run
an MCMC algorithm to obtain a chain $X^{(0)}=(X_{n}^{(0)})_{n\geq0}$
targeting $\pi^{(0)}$. Note that here the ``time'' index $n$
corresponds to
the number of iterations of the i-MCMC algorithm. We use the occupation
measure of the chain $X^{(0)}$ at time $n$ judiciously to design a second
MCMC algorithm to generate $X^{(1)}=(X_{n}^{(1)})_{n\geq0}$ at level 1
targeting $\pi^{(1)}$ which is typically more complex than $\pi^{(0)}$.
More precisely, the elementary transition $X_{n}^{(1)}\leadsto X_{n+1}^{(1)}$
of the chain $X^{(1)}$ at time $n$ depends on the occupation measure of
$%
( X_{0}^{(0)},X_{1}^{(0)},\ldots,X_{n}^{(0)} ) $.
Similarly we
use the empirical measure of $X^{(l-1)}$ at level $l-1$ to ``feed'' an MCMC
algorithm generating $X^{(l)}$ targeting $\pi^{(l)}$ at level $l$. These
i-MCMC samplers are SIMC in reference to the fact that the complete Markov
chain $\overline{X}_{n}^{m}:=(X_{n}^{(l)})_{0\leq l\leq m}$ associated with
a fixed series of $m$ levels evolves with elementary transitions
$\overline{X%
}_{n}^{m}\leadsto\overline{X}_{n+1}^{m}$ that depend on the occupation
measure of the whole system $\overline{X}_{p}^{m}$ from time $0$ up to
time $%
n$.

From the pure mathematical point of view, the convergence analysis of SIMC
is essentially based on the study of the stability properties of
sophisticated Markov chains with elementary transitions depending in a
nonlinear way on the occupation measure of the chains. Hence the theoretical
analysis of SIMC is much more involved than the one of traditional Markov
chains. It also differs significantly from interacting particle methods
developed in \cite{fk}. Besides the introduction of a new methodology, our
main contribution is a refined theoretical analysis based on measure-valued
processes and semigroup methods to analyze their asymptotic behavior as the
time index $n$ tends to infinity.

The rest of the paper is organized as follows:

The main notation used in this work are introduced in a brief
preliminary Section \ref{notat}. The i-MCMC methodology is detailed
formally in Section \ref{dmm}. The main results of the article are
presented in Section \ref{mrs}. Several examples of i-MCMC methods are
provided in Section \ref{exampl}. This section also provides a
discussion on how to combine interacting particle methods with i-MCMC
methods. Section \ref{timc} is concerned with the asymptotic behavior
of an abstract class of time inhomogeneous Markov chains. In Section
\ref{resl}, we present a preliminary resolvent analysis to estimate the
regularity properties of Poisson operator and invariant measure type
mappings. In Section \ref{lln}, we apply these results to study the law
of large numbers and the concentration properties of time inhomogeneous
Markov chains. In Section \ref{ladeff} we discuss the regularity
properties of a sequence of time averaged semigroups on distribution
flow state spaces. The asymptotic analysis of i-MCMC methods is
discussed in Section \ref{asymp}. The strong law of large numbers is
presented in Section~\ref{slln}. We also provide an
$\mathbb{L}_{r}$-mean error bound for the occupation measures of the
i-MCMC algorithms at each level $l$. In Section \ref{unifth}, we
discuss the long time behavior of these stochastic models in terms of
the exponential stability properties of a time averaged type semigroup
associated with the sequence of target measures. We prove a uniform
convergence theorem with respect to the level index $l$. The asymptotic
analysis of the occupation measures associated with the complete
self-interacting model on a fixed series of levels is discussed in
Section \ref{pathmodels}. The $\mathbb{L}_{r}$-mean error bounds and
the concentration analysis are presented, respectively, in
Sections \ref{lrmean} and in \ref{concentrat}. The final section,
Section \ref{FKS}, is concerned with contraction properties
of time averaged Feynman--Kac distribution flows.

\subsection{Notation and conventions}
\label{notat}

For the convenience of the reader we have collected some of
the main notation used in the article. We also recall some regularity
properties of integral operators used further in the article.

We denote, respectively, by $\mathcaligr{M}(E)$, $\mathcaligr{M}_{0}(E)$,
$\mathcaligr{P%
}(E)$ and $\mathcaligr{B}(E)$, the set of all finite signed measures on some
measurable space $(E,\mathcaligr{E})$, the convex subset of measures with null
mass, the set of all probability measures, and the Banach space of all
bounded and measurable functions $f$ on $E$. We equip $\mathcaligr{B}(E)$ with
the uniform norm $\Vert f\Vert={\sup_{x\in E}}|f(x)|$. We also denote
by $%
\mathcaligr{B}_{1}(E)\subset\mathcaligr{B}(E)$ the unit ball of functions
$f\in
\mathcaligr{B}(E)$ with $\Vert f\Vert\leq1$, and by $\mbox
{Osc}_{1}(E)$, the
convex set of $\mathcaligr{E}$-measurable functions $f$ with oscillations less
than one; that is,
\[
\mbox{osc}(f)=\sup{\{|f(x)-f(y)| ; x,y\in E\}}\leq1.
\]
We let $\mu(f)=\int\mu(dx) f(x)$ be the Lebesgue integral of a function
$f\in\mathcaligr{B}(E)$, with respect to a measure $\mu\in\mathcaligr{M}(E)$.
We slightly abuse the notation and sometimes denote by $\mu(A)=\mu(1_{A})$
the measure of a measurable subset $A\in\mathcaligr{E}$.

Let $M(x,dy)$ be a kernel from a measurable space $(E,\mathcaligr{E})$
into a
measurable space $(F,\mathcaligr{F})$ of the bounded integral operator $%
f\mapsto M(f)$ from $\mathcaligr{B}(F)$ into $\mathcaligr{B}(E)$ such that the
functions
\[
M(f)(x)=\int_{F} M(x,dy) f(y) \in\mathbb{R}
\]
are $\mathcaligr{E}$-measurable and bounded, for any $f\in\mathcaligr{B}(F)$.
Such a kernel also generates a dual operator $\mu\mapsto\mu M$ from $%
\mathcaligr{M}(E)$ into $\mathcaligr{M}(F)$ defined by $(\mu M)(f):=\mu(M(f))$.

We denote by $\Vert M\Vert:={\sup_{f\in\mathcaligr{B}_{1}(F)}}\Vert
M(f)\Vert$ the norm of the operator $f\mapsto M(f)$ and we equip the Banach
space $%
\mathcaligr{M}(E)$ with the corresponding total variation norm $\Vert\mu
\Vert
={\sup_{f\in\mathcaligr{B}_{1}(E)}}|\mu(f)|$. Using this slightly abusive
notation, we have
\[
\Vert M\Vert:={\sup_{x\in E}\sup_{f\in\mathcaligr{B}_{1}(F)}}|\delta
_{x}M(f)|%
={\sup_{x\in E}}\Vert\delta_{x}M\Vert,
\]
where $\delta_{x}$ stands for the Dirac measure at the point $x\in E$. We
recall that the norm of any kernel $M$ with null mass $M(1)=0$ satisfies
\[
\Vert M\Vert={\sup_{f\in\mathcaligr{B}_{1}(F)}}\Vert M(f)\Vert={2 \sup
_{f\in%
\mathrm{Osc}_{1}(F)}}\Vert M(f)\Vert.
\]
When $M$ has a constant mass, that is, $M(1) ( x )
=M(1) (
y ) $ for any $(x,y)\in E^{2}$, the operator $\mu\mapsto\mu M$
maps $%
\mathcaligr{M}_{0}(E)$ into $\mathcaligr{M}_{0}(F)$. In this situation, we
let $%
\beta(M)$ be the Dobrushin coefficient of a kernel $M$ defined by the
following formula:
\[
\beta(M):=\sup{\{\mathrm{osc}(M(f)) ; f\in\mathrm{Osc}_{1}(F)\}}.
\]
By construction, we have $M(f)/\beta(M)\in\mathrm{Osc}_{1}(E)$ as soon
as $\beta(M)\not=0$, so that
\begin{eqnarray*}
&&\Vert\mu M\Vert=2 {\sup_{f\in\mathrm{Osc}_{1}(F)}}|\mu
M(f)|\leq
{\beta(M) 2 \sup_{f\in\mathrm{Osc}_{1}(E)}}|\mu
(f)|
\\
&&\quad\Longrightarrow\quad
\Vert\mu M\Vert\leq\beta(M) \Vert\mu\Vert.
\end{eqnarray*}
Using the fact that $\Vert\delta_{x}-\delta_{y}\Vert=2$ for $x\neq y$
and
\begin{eqnarray*}
\beta(M) &=&{\sup_{f\in\mathrm{Osc}_{1}(F)}\sup_{(x,y)\in E^{2}}
}|(\delta_{x}M-\delta_{y}M)(f)|=\sup_{(x,y)\in E^{2}}{\frac{\Vert
\delta
_{x}M-\delta_{y}M\Vert}{\Vert\delta_{x}-\delta_{y}\Vert}} \\
&\leq&\sup_{\mu\in\mathcaligr{M}_{0}(E)}{\frac{\Vert\mu M\Vert
}{\Vert\mu
\Vert}}
\end{eqnarray*}
we prove that
\[
\beta(M)=\sup_{\mu\in\mathcaligr{M}_{0}(E)}\frac{\Vert\mu M\Vert
}{\Vert
\mu\Vert}={\frac{1}{2} \sup_{(x,y)\in E^{2}}}\Vert\delta
_{x}M-\delta
_{y}M\Vert
\]
is also the norm of the kernel
\[
\mu\in\mathcaligr{M}_{0}(E)\quad\mapsto\quad\mu M\in\mathcaligr{M}_{0}(F).
\]
That is, we have
\[
\beta(M)=\sup_{\mu\in\mathcaligr{M}_{0}(E)} ( {\Vert\mu M\Vert
}/{%
\Vert\mu\Vert} ) .
\]
More generally, for every kernel $K$ from a measurable space
$(E^{\prime},%
\mathcaligr{E}^{\prime})$ into an measurable space $(E,\mathcaligr{E})$, with
null mass $K(1)=0$, we have
\[
\Vert KM\Vert={\sup_{x\in E^{\prime}}}\Vert(\delta_{x}K)M\Vert
\leq
{\beta(M) \sup_{x\in E^{\prime}}}\Vert(\delta_{x}K)\Vert
\quad\Longrightarrow\quad
\Vert KM\Vert\leq\beta(M) \Vert K\Vert.
\]

Unless otherwise stated, we use the letter $C$ to denote a universal
constant whose value may vary from line to line. Finally, we shall use the
conventions $\sum_{\varnothing}=0$ and $\prod_{\varnothing}=1$.

\subsection{Interacting Markov chain Monte Carlo methods}
\label{dmm}

We describe here the \mbox{i-MCMC} methodology to numerically solve (\ref{phi}).
We consider a Markov transition $M^{(0)}$ from $S^{(0)}$ into itself
and a
collection of Markov transitions $M_{\mu}^{(l)}$ from $S^{(l)}$ into
itself, indexed by the parameter $l\geq0$ and the set of probability
measures $\mu\in\mathcaligr{P}(S^{(l-1)})$. We further assume that the
invariant measure of each operator $M_{\mu}^{(l)}$ is given by $\Phi
_{l}(\mu)$; that is, we have
\[
\forall l\geq0, \forall\mu\in\mathcaligr{P}\bigl(S^{(l-1)}\bigr)\qquad \Phi
_{l}(\mu)M_{\mu}^{(l)}=\Phi_{l}(\mu).
\]
For $l=0$, we use the convention $\Phi_{0}(\pi^{(-1)})=\pi^{(0)}$
and $%
M_{\mu}^{(0)}=M^{(0)}$. For every $l\leq m$, we denote by $\eta
^{(l)}\in
\mathcaligr{P}(S^{(l)})$ the image measure of a measure $\eta\in
\mathcaligr{P}%
(\prod_{0\leq l\leq m}S^{(l)})$ on the $l$th level space $S^{(l)}$.
We also
fix a sequence of probability measures $\nu_{k}$ on $S^{(k)}$, with
$k\geq0
$.

We let $X^{(0)}:=(X_{n}^{(0)})_{n\geq0}$ be a Markov chain on $S^{(0)}$
with initial distribution $\nu_{0}$ and Markov transitions $M^{(0)}$. For
every $k\geq1$, given a realization of the chain $%
X^{(k-1)}:=(X_{n}^{(k-1)})_{n\geq0}$, the $k$th level chain $X_{n}^{(k)}$
is a Markov chain with initial distribution $\nu_{k}$ and with random
Markov transitions\vspace*{-1pt} $M_{\eta_{n}^{(k-1)}}^{(k)}$ depending on the current
occupation measures $\eta_{n}^{(k-1)}$ of the chain at level $(k-1)$; that
is, we have
%
\begin{equation} \label{itdeff}
\mathbb{P}\bigl(X_{n+1}^{(k)}\in dx | X^{(k-1)}, X_{n}^{(k)}\bigr)=M_{\eta
_{n}^{(k-1)}}^{(k)}(X_{n}^{k},dx)
\end{equation}
with
\[
\eta_{n}^{(k-1)}:=\frac{1}{n+1}\sum_{p=0}^{n}\delta_{X_{p}^{(k-1)}}.
\]
The rationale behind this is that the $k$th level chain $X_{n}^{(k)}$
behaves asymptotically as a Markov chain with time homogeneous
transitions $%
M_{\pi^{(k-1)}}^{(k)}$ as long as $\eta_{n}^{(k-1)}$ is a good
approximation of $\pi^{(k-1)}$.

In the special case where $M_{\mu}^{(k)}(x^{k},\cdot)=\Phi
_{k}(\mu)$, the $k$th level chain $(X_{n}^{(k)})_{n\geq1}$ is a
collection of conditionally independent random variables with
distributions $%
(\Phi_{k}(\eta_{n-1}^{(k-1)}))_{n\geq1}$; that is, we have
%
\begin{equation} \label{imcmc0}\quad
\mathbb{P}\bigl(\bigl(X_{1}^{(k)},\ldots,X_{n}^{(k)}\bigr)\in
dx | X^{(k-1)}\bigr)=\prod_{p=1}^{n}\Phi_{k} \biggl( \frac{1}{p} \sum
_{0\leq
q<p} \delta_{X_{q}^{(k-1)}} \biggr) (dx_{p}),
\end{equation}
where $dx=d(x_{1},\ldots,x_{n})=dx_{1}\times\cdots\times dx_{n}$ stands
for an infinitesimal neighborhood of a generic path sequence $%
(x_{1},\ldots,x_{n})\in(S^{(k)})^{n}$.

We end this section with a SIMC interpretation of the stochastic algorithm
discussed above. We consider the product space
\[
E_{m}:=S^{(0)}\times\cdots\times S^{(m)}
\]
and we let $(K_{\eta}^{(m)})_{\eta\in\mathcaligr{P}(E_{m})}$ be the
collection of Markov transitions from $E_{m}$ into itself given by
%
\begin{equation}\label{defK}\quad
\forall x:=(x^{0},\ldots,x^{m})\in E_{m}\qquad K_{\eta
}^{(m)}(x,dy)=\prod_{0\leq l\leq
m}M_{\eta^{(l-1)}}^{(l)}(x^{l},dy^{l}),
\end{equation}
where $dy:=dy^{0}\times\cdots\times dy^{m}$ stands for for an
infinitesimal neighborhood of a generic point $y:=(y^{0},\ldots
,y^{m})\in
E_{m}$, and $\eta^{(l)}\in\mathcaligr{P}(S^{(l)})$ stands for the image
measure of a measure $\eta\in\mathcaligr{P}(E_{m})$ on the $l$th level space
$S^{(l)}$, with $m\geq l$. In other words, $\eta^{(l)}$ is the $l$th
marginal of the measure $\eta$. In this notation, we can readily check that
\[
\overline{X}_{n}^{m}:=\bigl(X_{n}^{(0)},\ldots,X_{n}^{(m)}\bigr)
\]
is an $E_{m}$-valued SIMC with elementary transitions defined by
%
\begin{equation} \label{pathocc}\qquad
\mathbb{P}(\overline{X}_{n+1}^{m}\in dy | \mathcaligr{F}_{n}^{\overline
{X}%
^{m}})=K_{\overline{\eta}_{n}^{[m]}}^{(m)}(\overline
{X}_{n}^{m},dy) \qquad
\mbox{with } \overline{\eta}_{n}^{[m]}=\frac{1}{n+1}%
\sum_{p=0}^{n}\delta_{\overline{X}_{n}^{m}},
\end{equation}
where $\mathcaligr{F}_{n}^{\overline{X}^{m}}$ stands for the filtration
generated by $\overline{X}^{m}$.

\subsection{Statement of some results}
\label{mrs}

We further\vspace*{1pt} assume that the mappings $\Phi_{l}\dvtx\mathcaligr{P}
(S^{(l-1)})\rightarrow\mathcaligr{P}(S^{(l)})$ satisfy the following
regularity condition for any $l\geq1$ and any pair of measures $(\mu
,\nu
)\in\mathcaligr{P}(S^{(l-1)})^{2}$
%
\begin{eqnarray}\label{lipcn}
&&\forall l\geq0, \forall f\in\mathcaligr{B}\bigl(S^{(l)}\bigr)
\nonumber\\[-8pt]\\[-8pt]
&&\qquad\vert%
[ \Phi_{l}(\mu)-\Phi_{l}(\nu) ] (f) \vert\leq
\int
\vert[ \mu-\nu] (g) \vert\Gamma_{l}(f,dg)\nonumber
\end{eqnarray}
for some kernel $\Gamma_{l}$ from $\mathcaligr{B}(S^{(l)})$ into
$\mathcaligr{B}%
(S^{(l-1)})$, with
\[
\int_{\mathcaligr{B}(S^{(l-1)})} \Gamma_{l}(f,dg) \Vert g\Vert\leq
\Lambda
_{l} \Vert f\Vert\quad\mbox{and}\quad\Lambda_{l}<\infty.
\]
We also suppose that there exist some integer $n_{l}\geq0$ and some
constant $c_{l}$ such that we have
%
\begin{equation} \label{condin}\hspace*{30pt}
\bigl\Vert M_{\mu}^{(l)}-M_{\nu}^{(l)}\bigr\Vert\leq c_{l} \Vert\mu-\nu
\Vert
\quad\mbox{and}\quad b_{l}(n_{l}):=\sup_{\mu\in\mathcaligr{P}(S^{(l-1)})}{\beta
\bigl(\bigl(M_{\mu
}^{(l)}\bigr)^{n_{l}}\bigr)}<1.
\end{equation}

This pair of abstract regularity conditions are rather standard. The first
one (\ref{lipcn}) is a natural Lipschitz property on the weakly continuous
integral mappings
\[
\forall f\in\mathcaligr{B}\bigl(S^{(l)}\bigr)\qquad \mu\in\mathcaligr{P}%
\bigl(S^{(l-1)}\bigr)\mapsto\Phi_{l}(\mu)(f)\in\mathbb{R}.
\]
Roughly speaking, this weak Lipschitz property simply expresses the fact
that $\Phi_{l}(\mu)(f)$ only depends on integrals of functions with
respect to the reference measure $\mu$. This condition is clearly satisfied
for linear Markov semigroups $\Phi_{l}(\mu)=\mu K_{l}$ associated with
some Markov transition $K_{l}$. We shall discuss this condition in the
context of nonlinear Feynman--Kac type semigroups (\ref{fksg}) in
Section \ref{sectFKmod}.

In the special case where $M_{\mu}^{(l)}(x^{l},\cdot)=\Phi
_{l}(\mu)$, the second condition (\ref{condin}) is trivially met for $
n_{l}=1$ with $b_{l}(n_{l})=0$. In this particular situation, the first
Lipschitz property of the mapping $\Phi_{l}(\mu)$ takes the following form:
\[
\Vert\Phi_{l}(\mu)-\Phi_{l}(\nu)\Vert\leq c_{l} \Vert\mu-\nu
\Vert.
\]
For more general models, condition (\ref{condin}) expresses the fact that
the Markov transitions $M_{\mu}^{(l)}$ are strongly continuous and they
satisfy Dobrushin's mixing condition, uniformly with respect to $\mu$. We
shall discuss this regularity condition in the context of
Metropolis--Hastings type algorithms (\ref{methasmod}) in Section \ref%
{methasmods}.

Under the conditions (\ref{condin}), for every $\eta\in\mathcaligr
{P}(E_{m})$%
, the invariant measure $\omega_{K_{\eta}^{(m)}}(\eta)\in\mathcaligr
{P}%
(E_{m})$ of $K_{\eta}^{(m)}$ defined in (\ref{defK}) is given by the tensor
product measure
%
\begin{equation} \label{invpath}
\omega_{K_{\eta}^{(m)}}(\eta)=\pi^{(0)}\otimes\Phi_{1}\bigl(\eta
^{(0)}\bigr)\otimes\cdots\otimes\Phi_{m}\bigl(\eta^{(m-1)}\bigr).
\end{equation}
We observe that the tensor product measure
%
\begin{equation}\label{eq:tensorproducttargetmeasure}
\overline{\pi}{}^{[m]}:=\pi^{(0)}\otimes\cdots\otimes\pi^{(m)}
\end{equation}
is a fixed point of the mapping $\omega_{K_{\eta}^{(m)}}\dvtx\eta\in
\mathcaligr{%
P}(E_{m})\rightarrow\omega_{K_{\eta}^{(m)}}(\eta)\in\mathcaligr{P}(E_{m})$.

Using this notation, our main results are basically as follows.
\begin{theorem}
\label{theointro} For any $r\geq1$, $m\geq1$, and any function $f\in
\mathcaligr{B}(E_{m})$ we have
\[
\sup_{n\geq1}\sqrt{n} \,\mathbb{E} \bigl( \bigl\vert\overline
{\eta}%
_{n}^{[m]}(f)-\overline{\pi}{}^{[m]}(f) \bigr\vert^{r} \bigr)
<\infty.
\]
Under some additional regularity conditions, we have the exponential
inequality
\[
\forall t>0\qquad \limsup_{n\rightarrow\infty}\frac{1}{n}\log
\mathbb{P}%
\bigl( \bigl\vert\bigl[ \overline{\eta}_{n}^{[m]}-\overline{\pi}{}^{[m]}%
\bigr] (f) \bigr\vert>t \bigr) <-\frac{t^{2}}{2\overline{\sigma
}_{m}^{2}}
\]
for some finite constant $\overline{\sigma}_{m}<\infty$ as well as the
following uniform convergence estimate:
\[
\sup_{k\geq0}\sup_{n\geq1}n^{\alpha/2} \mathbb{E} \bigl(
\bigl\vert\eta
_{n}^{(k)}(f_{k})-\pi^{(k)}(f_{k}) \bigr\vert^{r} \bigr) <\infty
\]
for some parameter $\alpha\in( 0,1 ] $ and for any
collection of
functions $(f_{k})_{k\geq0}\in\prod_{k\geq0}\mathcaligr{B}_{1}(S^{(k)})$.
\end{theorem}

We end this introduction with a series of comments and open research
questions.

First, the mean error bounds and the exponential estimates presented above
suggest the existence of Gaussian fluctuations of the occupation
measures $%
\overline{\eta}_{n}^{[m]}$ around their limiting value $\overline
{\pi}%
^{[m]}$, with a fluctuation rate $\sqrt{n}$. We have recently studied these
fluctuations in \cite{bercu2008,bercu20082}.

It might be surprising that the decays to equilibrium presented in
Theorem %
\ref{theointro} differ from the three types of decays exhibited in
\cite{dm1,dm2}. To understand the main differences between these classes of
interacting processes, we recall that the decay rate to equilibrium often
depends on the contraction coefficient of the invariant measure mapping
associated with a given self-interacting model. In our context, these
mappings are not necessarily contractive. Nevertheless, we shall see in
Section \ref{pathmodels} that the semigroup associated with these mappings
becomes essentially constant after a sufficiently large number of
iterations. In this respect, the self-interacting models discussed in the
present article are more regular than the ones analyzed in \cite{dm1,dm2}.

The uniform convergence estimate with respect to the number of levels
depends on the stability properties of a time averaged semigroup associated
with the mappings $\Phi_{l}$. The contraction properties of this new class
of nonlinear semigroups are studied in Section \ref{FKS} in the
context of
Feynman--Kac models. We show that the stability properties of the reference
Feynman--Kac semigroups can be transferred to study the associated time
averaged models. In more general situations this question remains open.

\section{Motivating applications}
\label{exampl}

\subsection{Feynman--Kac models}
\label{sectFKmod}

The main example of mappings $\Phi_{l}$ considered here are the Feynman--Kac
transformations given below:
%
\begin{eqnarray}\label{fksg}
&&\forall l\geq0, \forall(\mu,f)\in\bigl(\mathcaligr{P}\bigl(S^{(l)}\bigr)\times
\mathcaligr{B}\bigl(S^{(l+1)}\bigr)\bigr)\nonumber\\[-8pt]\\[-8pt]
&&\qquad\Phi_{l+1}(\mu)(f):={\mu
(G_{l}L_{l+1}(f))}/{%
\mu(G_{l})},\nonumber
\end{eqnarray}
where $G_{l}$ is a positive potential function on $S^{(l)}$, and $L_{l+1}$
stands for a Markov transition from $S^{(l)}$ into $S^{(l+1)}$. In this
situation, the solution of the measure-valued equation (\ref{phi}) is given
by the normalized Feynman--Kac distribution flow described below:
\[
\pi^{(l)}(f)={\gamma^{(l)}(f)}/{\gamma^{(l)}(1)}
\qquad\mbox{with }
\gamma^{(l)}(f):=\mathbb{E} \biggl( f(Y_{l}) \prod_{0\leq
k<l}G_{k}(Y_{k}) \biggr),
\]
where $(Y_{l})_{l\geq0}$ stands for a Markov chain taking values in the
state spaces $(S^{(l)})_{l\geq0}$, with initial distribution $\pi^{(0)}$
and Markov transitions $(L_{l})_{l\geq1}$. These probabilistic models arise
in a very wide variety of applications including nonlinear filtering and
rare event analysis as well the spectral analysis of Schroedinger type
operators and directed polymer analysis \cite{fk}. We also underline that
the unnormalized measures $\gamma^{(l)}$ are expressed in terms of
integrals on path spaces and we recall that $\gamma^{(l)}$ can be expressed
in terms of the sequence of measures $(\pi^{(k)})_{0\leq k<l}$ with the
following formulae:
%
\begin{equation}\label{muliteq}
\gamma^{(l)}(f)=\pi^{(l)}(f) \prod_{0\leq k<l}\pi^{(k)}(G_{k}).
\end{equation}
To check this assertion, we simply observe that
\[
\gamma^{(l)}(f)=\pi^{(l)}(f)\times\gamma^{(l)}(1)
\]
and we have the key multiplicative formula
%
\begin{eqnarray}\label{eq:normalizingconstants}
&&\gamma^{(l)}(1)=\gamma^{(l-1)}(G_{l-1})=\pi^{(l-1)}(G_{l-1})\times
\gamma
^{(l-1)}(1)\nonumber\\[-8pt]\\[-8pt]
&&\quad\Longrightarrow\quad\gamma^{(l)}(1)=\prod_{0\leq k<l}\pi
^{(k)}(G_{k}).\nonumber
\end{eqnarray}
Thus the i-MCMC methodology allows us to estimate the normalizing
constants $%
\gamma^{(l)}(1)$ by replacing the measures $\pi^{(k)}$ by their
approximations in (\ref{eq:normalizingconstants}). These models are quite
flexible. For instance, the reference Markov chain may represent the paths
from the origin up to the current time $l$ of an auxiliary chain $%
Y_{l}^{\prime}$ taking values in some state spaces $E_{l}^{\prime}$ with
some Markov transitions $(\widetilde{L}_{l})_{l\geq1}$ and potentials
$(%
\widetilde{G}_{l})_{l\geq1}$; that is, we have
%
\begin{equation} \label{pathmdels}
Y_{l}:=(Y_{0}^{\prime},\ldots,Y_{l}^{\prime})\in
S^{(l)}:=(E_{0}^{\prime
}\times\cdots\times E_{l}^{\prime})
\end{equation}
and
%
\begin{eqnarray}\label{eq:markovpotentialpathspaces}
L_{l} ( y_{l-1},d\overline{y}_{l} ) &=& \delta_{ (
y_{0}^{\prime
},\ldots,y_{l-1}^{\prime} ) } ( d ( \overline
{y}_{0}^{\prime
},\ldots,\overline{y}_{l-1}^{\prime} ) ) \widetilde
{L}_{l} (
y_{l-1}^{\prime},d\overline{y}_{l}^{\prime} ),\nonumber\\[-8pt]\\[-8pt]
G_{l} ( y_{l} )
&=& \widetilde{G}_{l} ( y_{l}^{\prime
} ).\nonumber
\end{eqnarray}

\subsection{Interacting Markov chain Monte Carlo methods for Feynman--Kac
models}\label{methasmods}

In the Feynman--Kac context and assuming we are working on path spaces
(\ref%
{pathmdels}), we can propose the following two i-MCMC algorithms to
approximate $\pi^{(l)}$. The first one simply consists of sampling directly
$X_{p}^{(k)}= ( X_{p}^{\prime(0)},X_{p}^{\prime
(1)},\ldots,X_{p}^{\prime
(k)} ) $ from the right-hand side product of the formula (\ref{imcmc0})
which takes here the following form:
\[
\Phi_{k} \biggl( \frac{1}{p} \sum_{0\leq q<p} \delta
_{X_{q}^{(k-1)}} \biggr)
\bigl( dx_{p}^{(k)} \bigr) =\sum_{0\leq q<p}\frac
{G_{k-1}(X_{q}^{(k-1)})}{%
\sum_{0\leq
m<p}G_{k-1}(X_{m}^{(k-1)})} L_{k}\bigl(X_{q}^{(k-1)},dx_{p}^{ (
k ) }\bigr),
\]
where $dx_{p}^{(k)}=dx_{p}^{\prime( 0 ) }\times\cdots
\times
dx_{p}^{\prime( k ) }$. We see that $X_{p}^{(k)}$ is sampled
according to two separate genetic type mechanisms. First, we randomly select
one state $X_{q}^{(k-1)}$ at level $(k-1)$ with a probability proportional
to its potential value $G_{k-1}(X_{q}^{(k-1)})$. Second, we randomly evolve
from this state according to the mutation transition~$L_{k}$. This i-MCMC
model can be interpreted as a spatial branching and interacting
process. In
this interpretation, the $k$th chain tends to duplicate individuals with
large potential values, at the expense of individuals\vspace*{1pt} with low potential
values. The selected offspring randomly evolve from the state space $%
S^{(k-1)}$ to the state space $S^{(k)}$ at the next level.

For the Feynman--Kac transformations (\ref{fksg}), we proved in \cite{fk}
that the condition (\ref{condin}) ensuring convergence of the
algorithm is
satisfied with $c_{l}=\beta(\widetilde{L}_{l})/\varepsilon_{l-1}(G)$ as soon
as the potential functions satisfy the following condition:\vspace*{8pt}

(G) \textit{For any $l\geq0$, the potential
functions $G_{l}$ are bounded above and bounded away from zero, so that}
\[
\varepsilon_{l}(G):=\inf_{x,y} \frac{G_{l}(x)}{G_{l}(y)}\in(0,1).
\]
\vspace*{0pt}

We can also propose the following alternative i-MCMC algorithm to
approximate $\pi^{(l)}$ which relies on using a transition kernel
$M_{\mu
}^{ ( l ) }$ different from $\Phi_{l} ( \mu)
$. We
introduce the following kernel from $S^{(l-1)}$ into $E_{l}^{\prime}$:
%
\begin{equation}\label{desintegration}
R_{l}((x_{0}^{\prime},\ldots,x_{l-1}^{\prime}),dx_{l}^{\prime})=%
\widetilde{L}_{l}(x_{l-1}^{\prime},dx_{l}^{\prime}) \widetilde{G}%
_{l-1}(x_{l-1}^{\prime}).
\end{equation}
In this scenario, it is sensible to propose to use for $M_{\mu}^{ (
l ) }$ in the i-MCMC algorithm the following Markov kernel on the
product space $S^{(l)}$ indexed by the set of measures $\mu\in
\mathcaligr{P}%
(S^{(l-1)})$
%
\begin{eqnarray}\label{methasmod}
M_{\mu}^{ ( l ) }(x,dy) &=& (\mu\otimes K_{l})(dy) \bigl(1\wedge r_{l}(x,y)\bigr)
\nonumber\\[-8pt]\\[-8pt]
&&{} + \biggl(1-\int_{S^{(l)}}\bigl(1\wedge r_{l}(x,z)\bigr) (\mu\otimes K_{l})(dz) \biggr)
\delta_{x}(dy),\nonumber
\end{eqnarray}
where $K_{l}$ is a Markov transition from $S^{(l-1)}$ into
$E_{l}^{\prime}$
and for every $(u,v)$ and $(w,z)\in(S^{(l-1)}\times E_{l}^{\prime})$
%
\begin{equation} \label{coondG}
r_{l}((u,v),(w,z)):=\frac{d ( K_{l}(u,\cdot)\otimes
R_{l}(w,%
\cdot) ) }{d ( R_{l}(u,\cdot)\otimes
K_{l}(w,%
\cdot) ) }(v,z),
\end{equation}
where we assume that
\[
K_{l}(u,\cdot)\otimes R_{l}(w,\cdot)\ll R_{l}(u,%
\cdot)\otimes K_{l}(w,\cdot).
\]
It can be checked that the kernel $M_{\mu}^{ ( l ) }$ is nothing
but a Metropolis--Hastings kernel of proposal distribution $\mu\otimes K_{l}$
and invariant distribution $\Phi_{l}(\mu)$.

We can also easily establish that for any measures $(\mu,\nu)\in
\mathcaligr{P%
}(S^{(l-1)})^{2}$%
\[
\bigl\Vert M_{\mu}^{(l)}-M_{\nu}^{(l)}\bigr\Vert\leq2 \Vert\mu-\nu\Vert
\]
so that the first condition on the left-hand side of (\ref{condin}) is
satisfied. Under the additional assumption that for any $ (
u,v )
\in(S^{(l-1)}\times E_{l}^{\prime})$
\[
\frac{dP_{l}(u,\cdot)}{dK_{l}(u,\cdot)} (
v )
\leq C_{l}
\]
it follows from \cite{mengersentweedie1996}, Theorem 2.1, that
\[
\beta\bigl(M_{\mu}^{ ( l ) }\bigr)\leq(1-C_{l}^{-1})
\]
from which we conclude that the second condition on the right-hand side
of (%
\ref{condin}) is met with $n_{l}=1$ and $b_{l}(n_{l})=(1-C_{l}^{-1})$.

\subsection{Interacting particle and Markov chain Monte Carlo methods}

As mentioned in the \hyperref[sec1]{Introduction}, in contrast to interacting particle
methods presented in Section \ref{ips}, we emphasize that the precision
parameter $n$ of i-MCMC models is not fixed but increases at every time
step. There exist several ways to combine an interacting particle method
with an i-MCMC method.

For instance,\vspace*{1pt} suppose we are given a realization of an interacting particle
algorithm $X^{(l)}=(X_{p}^{(l)})_{1\leq p\leq N}$ with a precision parameter
$N$. One natural way to initialize the i-MCMC model is to start with a
collection of initial random states $X_{0}^{(l)}$ sampled according to
the $%
N $-particle approximation measures
\[
\nu_{l}=\pi_{N}^{(l)}:=\frac{1}{N}\sum_{i=1}^{N}\delta_{X_{i}^{(l)}}.
\]
Another strategy is to use the $N$-particle approximation measures $\pi
_{N}^{(l)}$ in the evolution of the i-MCMC model. In other words, we
interpret the series of samples $X_{i}^{(l)}$, $1\leq i\leq N$, as the first
$N$ iterations of the i-MCMC model at level $l$. More formally, this
strategy simply substitutes the current occupation measure $\eta
_{n}^{(k-1)} $ of the chain at level $(k-1)$ in (\ref{itdeff}) by the
occupation measure $\eta_{n}^{(N,k-1)}$ of the whole sequence of random
variables at level $(k-1)$ defined by
\[
\eta_{n}^{(N,k-1)}=\frac{n+1}{N+n+1} \eta_{n}^{(k-1)}+\frac
{N}{N+n+1} \pi
_{N}^{(k-1)}.
\]

The convergence analysis of these two natural combinations of an interacting
particle method and i-MCMC method can be conducted easily using the
techniques developed in this article.

\section{Time inhomogeneous Markov chains}
\label{timc}

\subsection{Description of the models}

We consider a collection of Markov transitions $K_{\eta}$ on some
measurable space $(E,\mathcaligr{E})$ indexed by the set of probability
measures $\eta\in\mathcaligr{P}(F)$ on some possibly different measurable
space $(F,\mathcaligr{F})$. We further assume that for any pair of
measures $%
(\eta,\mu)\in\mathcaligr{P}(F)^{2}$ and some integer $n_{0}\geq0$ we have
%
\begin{equation} \label{c1}
\Vert K_{\eta}-K_{\mu}\Vert\leq c \Vert\eta-\mu\Vert
\quad\mbox{and}\quad b(n_{0}):=\sup_{\eta\in\mathcaligr{P}(E)}{\beta(K_{\eta
}^{n_{0}})}%
<1.
\end{equation}
We associate with the collection of transitions $K_{\eta}$ an $E$-valued
inhomogeneous random process $X_{n}$ with elementary transitions
defined by
\[
\mathbb{P}(X_{n+1}\in dx | X_{0},\ldots,X_{n})=K_{\mu_{n}}(X_{n},dx),
\]
where $\mu_{n}$ is a sequence of possibly random distributions on $F$ that
only depends on the random sequence $(X_{0},\ldots,X_{n})$. More precisely,
$\mu_{n}$ is a measurable random variable with respect to the
$\sigma$-field generated by the random states $X_{p}$ from the origin $p=0$, up to
the current time horizon $p=n$. We further assume that the variations
of the
flow $\mu_{n}$ are controlled by some sequence of random variables $%
\varepsilon(n)$ in the sense that
\[
\forall n\geq0\qquad \Vert\mu_{n+1}-\mu_{n}\Vert\leq\varepsilon(n).
\]
We let $\overline{\varepsilon}(n)$ be the mean variation of the distribution
flow $(\mu_{p})_{0\leq p\leq n}$; that is, we have
\[
\overline{\varepsilon}(n):=\frac{1}{n+1}\sum_{p=0}^{n}\varepsilon(p).
\]
For SIMC, we have $F=E$ and the measure $\mu_{n}$ coincides with the
occupation measures of the chain up to the current time $n$. In this
particular situation, we have
%
\begin{equation} \label{simc}
\mu_{n}=\eta_{n}:=\frac{1}{n+1} \sum_{p=0}^{n} \delta
_{X_{p}} \quad\Longrightarrow\quad\varepsilon(n)\leq\frac{2}{n+2}.
\end{equation}
This implies that
\[
\overline{\varepsilon}(n)\leq\frac{2}{n+1} \log{(n+2)}.
\]
Under assumption (\ref{c1}), every elementary transition $K_{\mu_{n}}(x,dy)$
admits an invariant measure
\[
\omega(\mu_{n})K_{\mu_{n}}=\omega(\mu_{n})\in\mathcaligr{P}(E).
\]
For sufficiently small variations $\varepsilon(n)$ of the distribution
flow $%
\mu_{n}$, we expect that the occupation measures $\eta_{n}$ have the same
asymptotic behavior as the mean values $\overline{\omega}_{n}(\mu)$
of the
instantaneous invariant measures $\omega(\mu_{p})$ from time $p=0$ up to
the current time $p=n$. That is, for large values of the time horizon $n$,
we have in some sense
%
\begin{equation} \label{averageinvariantdistributions}
\eta_{n}\simeq\overline{\omega}_{n}(\mu):=\frac{1}{n+1}%
\sum_{p=0}^{n}\omega(\mu_{p}).
\end{equation}

\subsection{A resolvent analysis}
\label{resl}

We recall that assumption (\ref{c1}) ensures that
$K_{\eta}$
has a unique invariant measure for any $\eta\in\mathcaligr{P}(F)$
\[
\omega(\eta)K_{\eta}=\omega(\eta)\in\mathcaligr{P}(E)
\]
and the pair of sums given by
%
\begin{equation}\label{ref2}
\alpha(\eta):=\sum_{n\geq0}\beta(K_{\eta}^{n})\in[
1,\infty
) \quad\mbox{and}\quad \sum_{n\geq0} [ K_{\eta}^{n}-\omega(\eta) ] (f)
\end{equation}
are absolutely convergent for any $f\in\mathcaligr{B}(E)$. The main
simplification of these conditions comes from the fact that the resolvent
operator
\[
P_{\eta} \dvtx f\in\mathcaligr{B}(E)\quad\rightarrow\quad P_{\eta}(f):=\sum_{n\geq
0} [
K_{\eta}^{n}-\omega(\eta) ] (f)\in\mathcaligr{B}(E)
\]
is a well-defined solution of the Poisson equation
\[
\cases{
(K_{\eta}-\mathrm{Id})P_{\eta} = \bigl(\omega(\eta)-\mathrm{Id}\bigr), \cr
\omega(\eta)P_{\eta} = 0.}
\]
The reader should not be misled by the notation $P_{\eta}$. In this
context, $P_{\eta}$ is not a Markov transition kernel. We have used the
letter $P$ in reference to the solution of the Poisson equation.
\begin{prop}
\label{lem1} For any $\eta\in\mathcaligr{P }(F)$, $P_{\eta}$ is a bounded
integral operator on $\mathcaligr{B }(E)$ and we have
\[
(\Vert P_{\eta}\Vert/2)\vee\beta(P_{\eta})\leq\alpha(\eta)\leq
\frac{n_{0}}{%
1-\beta(K^{n_{0}}_{\eta})}.
\]
\end{prop}
\begin{pf}
The fact that $\beta( P_{\eta} ) \leq\alpha(\eta)$ is readily
deduced from the following decomposition:
\[
P_{\eta}(f)(x)-P_{\eta}(f)(y):=\sum_{n\geq0} [ K_{\eta
}^{n}(f)(x)-K_{%
\eta}^{n}(f)(y) ].
\]
Indeed, using this decomposition we find that $\mathrm{osc}%
(P_{\eta}(f))\leq\sum_{n\geq0}\operatorname{osc}(K_{\eta}^{n}(f))$. Recalling
that $\mathrm{osc}(K_{\eta}^{n}(f))\leq\beta(K_{\eta}^{n}) \operatorname
{osc}%
(f)$, we conclude that
\[
\operatorname{osc}(P_{\eta}(f))\leq\biggl[ \sum_{n\geq0}\beta(K_{\eta
}^{n}) \biggr]
\operatorname{osc}(f)\quad\Rightarrow\quad\beta(P_{\eta})\leq\sum_{n\geq
0}\beta(K_{\eta}^{n}).
\]
In much the same way, we use the fact that
\[
P_{\eta}(f)(x)=\sum_{n\geq0}\int [ K_{\eta
}^{n}(f)(x)-K_{\eta}^{n}(f)(y)%
] \omega(\eta)(dy)
\]
to check that
\[
\Vert P_{\eta}(f)\Vert\leq\sum_{n\geq0}\operatorname{osc}(K_{\eta}^{n}(f))
\]
and
\[
\Vert P_{\eta}(f)\Vert\leq\biggl[ \sum_{n\geq0}\beta(K_{\eta
}^{n}) \biggr] %
\operatorname{osc}(f)\quad\Rightarrow\quad\Vert
P_{\eta}\Vert\leq2 \sum_{n\geq0}\beta(K_{\eta}^{n}).
\]
To prove that $\alpha(\eta)\leq\frac{n_{0}}{1-\beta(K_{\eta
}^{n_{0}})}$, we
use the decomposition
\[
\alpha(\eta):=\sum_{n\geq0}\beta(K_{\eta}^{n})=\sum_{p\geq
1} \sum
_{n=(p-1)n_{0}}^{pn_{0}-1}\beta(K_{\eta}^{n})=\sum_{p\geq1} %
\sum_{r=0}^{n_{0}-1}\beta\bigl(K_{\eta}^{(p-1)n_{0}+r}\bigr).
\]
Since we have
\[
\beta\bigl(K_{\eta}^{(p-1)n_{0}+r}\bigr)\leq\beta\bigl(K_{\eta
}^{(p-1)n_{0}}\bigr) \beta(K_{\eta
}^{r})\leq\beta(K_{\eta}^{n_{0}})^{(p-1)} \beta(K_{\eta}^{r})\leq
\beta
(K_{\eta}^{n_{0}})^{(p-1)}
\]
we conclude that $\alpha(\eta) \leq n_{0} \sum_{p\geq0} \beta(K_{%
\eta}^{n_{0}})^{p}=\frac{n_{0}}{1-\beta(K_{\eta}^{n_{0}})}$. The
end of the
proof of the proposition is now complete.
\end{pf}
\begin{prop}
\label{lem2} For any pair of measures $(\eta,\mu)\in\mathcaligr{P
}(F)^{2}$, we
have
%
\begin{equation} \label{e1}
\|\omega(\eta)-\omega(\mu)\|\leq\delta_{n_{0}}(\eta,\mu) \|
\eta-\mu\|
\end{equation}
and
\[
\Vert P_{\mu}-P_{\eta}\Vert\leq\alpha(\eta) [ 2c \alpha
(\mu
)+\delta_{n_{0}}(\eta,\mu) ] \|\eta-\mu\|
\]
for some finite constant $\delta_{n_{0}}(\eta,\mu)$ such that
%
\begin{equation} \label{ref3}
\delta_{n_{0}}(\eta,\mu)\leq\frac{cn_{0}}{1-(\beta
(K^{n_{0}}_{\eta})\wedge%
\beta(K^{n_{0}}_{\mu}))}.
\end{equation}
\end{prop}
\begin{pf}
The proof of the first assertion is based on the following decomposition:
\[
\omega(\eta)-\omega(\mu)=\omega(\eta)(K_{\eta}^{n_{0}}-K_{\mu}^{n_{0}})+
[ \omega(\eta)-\omega(\mu) ] K_{\mu}^{n_{0}}.
\]
Using the fact that
\[
\Vert[ \omega(\eta)-\omega(\mu) ] K_{\mu
}^{n_{0}}\Vert\leq
\beta(K_{\mu}^{n_{0}}) \Vert\omega(\eta)-\omega(\mu)\Vert
\]
we find that
%
\begin{equation}\label{ref1}
\Vert\omega(\eta)-\omega(\mu)\Vert\leq\frac{1}{1-(\beta(K_{\mu
}^{n_{0}})%
\wedge\beta(K_{\eta}^{n_{0}}))} \Vert\omega(\eta)(K_{\eta
}^{n_{0}}-K_{\mu
}^{n_{0}})\Vert.
\end{equation}
On the other hand, we have
\[
\Vert\omega(\eta)(K_{\eta}^{n_{0}}-K_{\mu}^{n_{0}})\Vert\leq
\Vert K_{\eta
}^{n_{0}}-K_{\mu}^{n_{0}}\Vert\Vert\omega(\eta)\Vert=\Vert
K_{\eta}^{n_{0}}-K_{\mu}^{n_{0}}\Vert.
\]
Using the decomposition
\[
K_{\eta}^{n_{0}}-K_{\mu}^{n_{0}}=\sum_{p=0}^{n_{0}-1} K_{\mu
}^{p}(K_{\eta
}-K_{\mu})K_{\eta}^{n_{0}-(p+1)}
\]
we find that
\[
\Vert K_{\eta}^{n_{0}}-K_{\mu}^{n_{0}}\Vert\leq\sum
_{p=0}^{n_{0}-1} \bigl\Vert
K_{\mu}^{p}(K_{\eta}-K_{\mu})K_{\eta}^{n_{0}-(p+1)}\bigr\Vert.
\]
For any $0\leq p\leq n_{0}$ we have
\begin{eqnarray*}
\bigl\Vert K_{\mu}^{p}(K_{\eta}-K_{\mu})K_{\eta}^{n_{0}-(p+1)}\bigr\Vert
&\leq&
\Vert K_{\mu}^{p}\Vert\Vert K_{\eta}-K_{\mu}\Vert\bigl\Vert
K_{\eta}^{n_{0}-(p+1)}\bigr\Vert\\
&\leq&
\Vert K_{\eta}-K_{\mu}\Vert\leq c \Vert\eta-\mu\Vert
\end{eqnarray*}
from which we conclude that
\[
\Vert K_{\eta}^{n_{0}}-K_{\mu}^{n_{0}}\Vert\leq cn_{0} \Vert\eta
-\mu\Vert\quad\Longrightarrow\quad\Vert\omega(\eta)(K_{\eta}^{n_{0}}-K_{\mu
}^{n_{0}})%
\Vert\leq cn_{0} \Vert\eta-\mu\Vert.
\]
The proof of (\ref{e1}) is now a direct consequence of (\ref{ref1}).

The proof of the second assertion is based on the following decomposition:
\[
P_{\eta}-P_{\mu}=P_{\mu}(K_{\eta}-K_{\mu})P_{\eta}+ [ \omega
(\mu
)-\omega(\eta) ] P_{\eta}.
\]
To check this formula, we first use the fact that $K_{\mu}P_{\mu
}=P_{\mu
}K_{\mu}$ to prove that
\[
P_{\mu}(K_{\mu}-\mathrm{Id})=(K_{\mu}-\mathrm{Id})P_{\mu}=\bigl(\omega(\mu)-\mathrm{Id}\bigr).
\]
This yields
\[
P_{\mu}(K_{\mu}-\mathrm{Id})P_{\eta}=\bigl(\omega(\mu)-\mathrm{Id}\bigr)P_{\eta}.
\]
Using the Poisson equation and using the fact that $P_{\mu}(1)=0$ we also
have the decomposition
\[
P_{\mu}(K_{\eta}-\mathrm{Id})P_{\eta}=P_{\mu}\bigl(\omega(\eta)-\mathrm{Id}\bigr)=-P_{\mu}.
\]
Combining these two formulae, we conclude that
\[
P_{\mu}(K_{\eta}-K_{\mu})P_{\eta}= [ P_{\eta}-P_{\mu}
] -%
[ \omega(\mu)-\omega(\eta) ] P_{\eta}.
\]
It follows that
\[
\Vert P_{\eta}-P_{\mu}\Vert\leq\Vert P_{\mu}(K_{\eta}-K_{\mu
})P_{\eta
}\Vert+\Vert[ \omega(\mu)-\omega(\eta) ] P_{\eta
}\Vert.
\]
The term on the right-hand side is easily estimated. Indeed, under our
assumptions we readily find that
\begin{eqnarray*}
\Vert[ \omega(\mu)-\omega(\eta) ] P_{\eta}\Vert
&\leq&
\beta(P_{\eta}) \Vert\omega(\eta)-\omega(\mu)\Vert\\
&\leq&
\alpha(\eta) \Vert\omega(\eta)-\omega(\mu)\Vert\leq
\alpha
(\eta) \delta_{n_{0}}(\eta,\mu) \Vert\eta-\mu\Vert.
\end{eqnarray*}
On the other hand, we have
\[
\Vert P_{\mu}(K_{\eta}-K_{\mu})P_{\eta}\Vert\leq\beta(P_{\eta
}) \Vert P_{\mu}(K_{\eta}-K_{\mu})\Vert\leq\beta(P_{\eta
}) \Vert
P_{\mu}\Vert\Vert K_{\eta}-K_{\mu}\Vert
\]
from which we conclude that
\[
\Vert P_{\mu}(K_{\eta}-K_{\mu})P_{\eta}\Vert\leq2c \alpha(\mu
)\alpha(\eta) \Vert\eta-\mu\Vert.
\]
The end of the proof is now clear.
\end{pf}

\subsection{$\mathbb{L}_{r}$-inequalities and concentration analysis}
\label{lln}

First, we examine some of the consequences of the pair of regularity
conditions presented in (\ref{c1}). The second condition ensures that the
functions $\alpha(\eta)$ and $\delta_{n_{0}}(\eta,\mu)$ introduced
in (\ref%
{ref2}) and (\ref{ref3}) are uniformly bounded; that is, we have
%
\begin{equation} \label{an0}
1\leq a(n_{0}):=\sup_{\eta\in\mathcaligr{P}(F)}{\alpha(\eta)}\leq
\frac{n_{0}}{%
1-b(n_{0})}
\end{equation}
and
%
\begin{equation}\label{dn0}
d(n_{0}):=\sup_{(\eta,\mu)\in\mathcaligr{P}(F)^{2}}{\delta
_{n_{0}}(\eta,\mu)}%
\leq\frac{cn_{0}}{1-b(n_{0})}<\infty.
\end{equation}
We recall that $\overline{\omega}_{n}(\mu)$ is defined in (\ref%
{averageinvariantdistributions}). We are now in a position to state and
prove the main result of this section.
\begin{theorem}
\label{theointer} For any $n\geq0$, $f\in\mathcaligr{B}_{1}(E)$ and
$r\geq1$ we
have the estimate
\[
\mathbb{E} \bigl( \vert[ \eta_{n}-\overline{\omega
}_{n}(\mu) ]
(f) \vert^{r} \bigr) ^{{1}/{r}}\leq e(r) \biggl( \frac
{n_{0}}{%
1-b(n_{0})} \biggr) ^{2} \biggl[ \frac{1}{\sqrt{n+1}}+c \mathbb
{E}(\overline{%
\varepsilon}(n)^{r})^{{1}/{r}} \biggr]
\]
for some finite constant $e(r)<\infty$ whose value only depends on the
parameter $r$. In addition, for any $\delta\in( 0,1 ) $
and any
time horizon $n\geq1$, the probability that
\begin{eqnarray*}
&& \vert[ \eta_{n}-\overline{\omega}_{n}(\mu) ]
(f) \vert
\\
&&\qquad\leq\frac{n_{0}}{1-b(n_{0})} \Biggl[ \sqrt{\frac{2\log{(2/\delta
)}}{n+1}}%
+(1+c) \biggl( \frac{4n_{0}}{1-b(n_{0})} \biggr) \biggl[ \overline
{\varepsilon}%
(n)\vee\frac{1}{n+1} \biggr] \Biggr]%
\end{eqnarray*}
is greater than $(1-\delta)$ [where $c$ is the constant introduced in
(\ref{c1})].
\end{theorem}
\begin{cor}
\label{corinter} For the SIMC associated with the occupation measure
distribution flow (\ref{simc}), we have for any $n\geq0$, $f\in
\mathcaligr{B}%
_{1}(E)$ and any $r\geq1$
\[
\sqrt{n+1} \mathbb{E} \bigl( \vert[ \eta_{n}-\overline
{\omega}%
_{n}(\mu) ] (f) \vert^{r} \bigr) ^{{1}/{r}}\leq
e(r) (1+c) \biggl( \frac{n_{0}}{1-b(n_{0})} \biggr) ^{2}
\]
for some finite constant $e(r)<\infty$ whose value only depends on the
parameter $r$. In addition, for any $\delta\in( 0,1 ) $
and any
time horizon $n\geq1$, the probability that
\[
\vert[ \eta_{n}-\overline{\omega}_{n}(\mu) ]
(f) \vert\leq\biggl( \frac{2n_{0}}{1-b(n_{0})} \biggr)
^{2} \sqrt{\frac{%
2}{n+1}} \bigl[ \sqrt{\log{(2/\delta)}}+2(1+c) \bigr]
\]
is greater than $(1-\delta)$.
\end{cor}
\begin{pf*}{Proof of Theorem \protect\ref{theointer}}
First, we examine some consequences of the regularity conditions
presented in (\ref{c1}) on the resolvent function $P_{\eta}$ introduced
in~(\ref{ref2}). Using Propositions \ref{lem1} and \ref{lem2} we find
the following uniform estimates:
\[
\sup_{\eta\in\mathcaligr{P}(F)} \bigl( (\Vert P_{\eta}\Vert/2)\vee
\beta
(P_{\eta}) \bigr) \leq\frac{n_{0}}{1-b(n_{0})}
\]
and
%
\begin{equation}\label{pmua}
\Vert P_{\mu}-P_{\eta}\Vert\leq3c \biggl( \frac
{n_{0}}{1-b(n_{0})} \biggr)
^{2} \Vert\mu-\eta\Vert.
\end{equation}
In addition, using Proposition \ref{lem2} again we find that the invariant
measure mapping $\omega$ is uniform Lipschitz in the sense that
\[
\Vert\omega(\eta)-\omega(\mu)\Vert\leq\frac{cn_{0}}{1-b(n_{0})}
\Vert\eta-\mu\Vert.
\]
For any $n\geq0$ and any function $f\in\mathcaligr{B}_{1}(E)$, we set
\[
I_{n}(f):=(n+1) [ \eta_{n}-\overline{\omega}_{n}(\mu) ]
(f)=\sum_{p=0}^{n} [ f(X_{p})-\omega(\mu_{p})(f) ].
\]
Using the Poisson equation, we have
\[
[ \mathrm{Id}-\omega(\mu_{p}) ] =(\mathrm{Id}-K_{\mu_{p}})P_{\mu_{p}}.
\]
From this formula, we find the decomposition
%
\begin{eqnarray} \label{deco}
&&[ f(X_{p})-\omega(\mu_{p})(f) ]\nonumber\\
&&\qquad= P_{\mu_{p}}(f)(X_{p})-K_{\mu
_{p}}(P_{\mu_{p}}(f))(X_{p}) \\
&&\qquad= [ P_{\mu_{p}}(f)(X_{p})-P_{\mu_{p}}(f)(X_{p+1}) ]
+\Delta
M_{p+1}(f)\nonumber
\end{eqnarray}
with the increments
\[
\Delta M_{p+1}(f):= [ P_{\mu_{p}}(f)(X_{p+1})-K_{\mu_{p}}(P_{\mu
_{p}}(f))(X_{p}) ]
\]
of the martingale $M_{n+1}(f)$ defined by
\[
M_{n+1}(f):=\sum_{p=1}^{n+1}\Delta M_{p}(f)=\sum_{p=1}^{n+1} [
P_{\mu
_{p-1}}(f)(X_{p})-K_{\mu_{p-1}}(P_{\mu_{p-1}}(f))(X_{p-1}) ].
\]
For $n=0$, we set $M_{0}(f)=0$. The first term in the right-hand side
of (%
\ref{deco}) can also be rewritten in the following form:
\begin{eqnarray*}
&&P_{\mu_{p}}(f)(X_{p})-P_{\mu_{p}}(f)(X_{p+1}) \\
&&\qquad= [ P_{\mu_{p}}(f)(X_{p})-P_{\mu_{p+1}}(f)(X_{p+1}) ]
\\
&&\qquad\quad{} + [
P_{\mu_{p+1}}(f)(X_{p+1})-P_{\mu_{p}}(f)(X_{p+1}) ].
\end{eqnarray*}
This yields the decomposition
\begin{eqnarray*}
&&\sum_{p=0}^{n} [ P_{\mu_{p}}(f)(X_{p})-P_{\mu
_{p}}(f)(X_{p+1}) ] \\
&&\qquad=
[ P_{\mu_{0}}(f)(X_{0})-P_{\mu_{n+1}}(f)(X_{n+1}) ] +L_{n+1}(f)
\end{eqnarray*}
with the random sequence
\[
L_{n+1}(f):=\sum_{p=0}^{n} [ P_{\mu_{p+1}}-P_{\mu_{p}} ]
(f)(X_{p+1}).
\]
In summary, we have established the following decomposition:
\[
I_{n}(f)=M_{n+1}(f)+L_{n+1}(f)+ [ P_{\mu_{0}}(f)(X_{0})-P_{\mu
_{n+1}}(f)(X_{n+1}) ].
\]
We estimate each term separately. First, using (\ref{pmua}) we prove that
\[
\vert P_{\mu_{0}}(f)(X_{0})-P_{\mu_{n+1}}(f)(X_{n+1})
\vert\leq
\Vert P_{\mu_{0}}\Vert+\Vert P_{\mu_{n+1}}\Vert\leq\frac{4n_{0}}{%
1-b(n_{0})}.
\]
In much the same way, using (\ref{pmua}) we obtain
\begin{eqnarray*}
\Vert L_{n+1}\Vert
&\leq&
\sum_{p=0}^{n}\Vert P_{\mu_{p+1}}-P_{\mu
_{p}}\Vert
\leq3c \biggl( \frac{n_{0}}{1-b(n_{0})} \biggr) ^{2} \sum
_{p=0}^{n}\Vert\mu
_{p+1}-\mu_{p}\Vert\\
&=&
3c (n+1) \biggl( \frac{n_{0}}{1-b(n_{0})} \biggr) ^{2} \overline
{\varepsilon}(n).
\end{eqnarray*}
From these two estimates, we conclude that
%
\begin{equation} \label{inter}\quad
\vert I_{n}(f) \vert\leq\vert M_{n+1}(f)
\vert
+3c (n+1) \biggl( \frac{n_{0}}{1-b(n_{0})} \biggr) ^{2} \overline
{\varepsilon}%
(n)+\frac{4n_{0}}{1-b(n_{0})}.
\end{equation}
To estimate the martingale term, we recall that the unpredictable quadratic
variation process $ [ M(f),M(f) ] _{n}$ of the martingale $M_{n}(f)$
is the cumulated sum of the square of its increments from the origin up to
the current time; that is, we have
\[
[ M(f),M(f) ] _{n}:=\sum_{p=1}^{n} (\Delta M_{p}(f))^{2}.
\]
The main simplification of our regularity conditions comes from the fact
that the increments $|\Delta M_{p}(f)|$ are uniformly bounded. More
precisely, we have the almost sure estimates
\begin{eqnarray*}
\vert\Delta M_{p+1}(f) \vert
&=&
\vert P_{\mu_{p}}(f)(X_{p+1})-K_{\mu_{p}}(P_{\mu_{p}}(f))(X_{p}) \vert\\
&=&
\biggl\vert\int [ P_{\mu_{p}}(f)(X_{p+1})-P_{\mu_{p}}(f)(x) ]
K_{\mu_{p}}(X_{p},dx) \biggr\vert\\
&\leq&
\int \vert P_{\mu_{p}}(f)(X_{p+1})-P_{\mu
_{p}}(f)(x) \vert K_{\mu_{p}}(X_{p},dx)
\end{eqnarray*}
from which we conclude that
\[
\vert\Delta M_{p+1}(f) \vert\leq\operatorname{osc}(P_{\mu
_{p}}(f))\leq\beta(P_{\mu_{p}})\leq\frac{n_{0}}{1-b(n_{0})}.
\]
By definition of the quadratic variation process $ [
M(f),M(f) ]
_{n}$, this implies that
\[
[ M(f),M(f) ] _{n}\leq\biggl( \frac
{n_{0}}{1-b(n_{0})} \biggr)
^{2} n.
\]
The end of the proof is now a direct consequence of the
Burkholder--Davis--Gundy inequality for martingales. For any $r\geq1$, there
exists some finite constant $e(r)$ whose value only depends on $r$, and such
that for any $n$
\[
\mathbb{E} \Bigl( {\max_{1\leq p\leq n}} \vert M_{p}(f)
\vert^{r}%
\Bigr) ^{1/r} \leq e(r) \mathbb{E} ( [
M(f),M(f) ]
_{n}^{r/2} ) ^{1/r}\leq e(r) \frac
{n_{0}}{1-b(n_{0})} %
\sqrt{n}.
\]
Combining this estimate with (\ref{inter}), we find that
\[
\mathbb{E} ( \vert I_{n}(f) \vert^{r} )
^{1/r}%
 \leq e(r) \biggl( \frac{n_{0}}{1-b(n_{0})} \biggr) ^{2} \bigl[
\sqrt{(n+1)}%
+c (n+1) \mathbb{E}(\overline{\varepsilon}(n)^{r})^{1/r} \bigr]
\]
with again some finite constant $e(r)$ whose values may vary from line to
line, but only depends on $r$. Recalling the definition of $I_{n}(f)$, we
conclude that
\[
\mathbb{E} \bigl( \vert[ \eta_{n}-\overline{\omega
}_{n}(\mu)%
] (f) \vert^{r} \bigr) ^{1/r}\leq e(r) \biggl(
\frac{n_{0}%
}{1-b(n_{0})} \biggr) ^{2} \biggl[ \frac{1}{\sqrt{(n+1)}}+c \mathbb
{E}(%
\overline{\varepsilon}(n)^{r})^{1/r} \biggr].
\]
This ends the proof of the first assertion. To prove the concentration
estimates, we use the fact that
\[
\vert[ \eta_{n}-\overline{\omega}_{n}(\mu) ]
(f) \vert\leq\frac{ \vert M_{n+1}(f) \vert
}{n+1}+\frac{n_{0}%
}{1-b(n_{0})} \biggl[ \frac{3c n_{0}}{1-b(n_{0})} \overline{\varepsilon
}(n)+%
\frac{4}{n+1} \biggr]
\]
from which we deduce the rather crude upper bound
%
\begin{eqnarray}\label{chernov}
&&\vert[ \eta_{n}-\overline{\omega}_{n}(\mu) ]
(f) \vert
\nonumber\\[-8pt]\\[-8pt]
&&\qquad\leq\frac{ \vert M_{n+1}(f) \vert}{n+1}%
+(1+c) \biggl( \frac{2n_{0}}{1-b(n_{0})} \biggr) ^{2} \biggl[ \overline
{\varepsilon
}(n)\vee\frac{1}{n+1} \biggr].\nonumber
\end{eqnarray}
The Chernov--Hoeffding exponential inequality states that for every
martingale $M_{n}$ with $M_{0}=0$ and uniformly bounded increments $%
{\sup_{n}}|\Delta M_{n}|\leq a$, we have
\[
\mathbb{P}( \vert M_{n} \vert\geq tn)\leq2 e^{-{nt^{2}}/{2a^{2}}}.
\]
In our context, we have proved that ${\sup_{n}}|\Delta M_{n}(f)|\leq
n_{0}/(1-b(n_{0}))$, from which we conclude that
\begin{eqnarray*}
&&\mathbb{P} \biggl( \vert[ \eta_{n}-\overline{\omega
}_{n}(\mu)%
] (f) \vert>t+(1+c) \biggl( \frac
{2n_{0}}{1-b(n_{0})} \biggr) ^{2}%
\biggl[ \overline{\varepsilon}(n)\vee\frac{1}{n+1} \biggr] \biggr) \\
&&\qquad
\leq2 \exp \biggl( -(n+1) \frac{t^{2}}{2} \biggl( \frac
{1-b(n_{0})}{n_{0}}%
\biggr) ^{2} \biggr).
\end{eqnarray*}
We conclude the proof of the theorem by choosing $t=\frac
{n_{0}}{1-b(n_{0})}%
\sqrt{\frac{2\log{(2/\delta)}}{n+1}}$.
\end{pf*}

\section{Distribution flows models}
\label{ladeff}

In this section, we have collected the definition of a series
of semigroups on distribution flow spaces. We also take the opportunity to
describe some of their regularity properties we shall use in the further
developments of the article.

We equip the sets of distribution flows $\mathcaligr{P}(S^{(l)})^{\mathbb{N}}$
with the uniform total variation distance defined by
\[
\forall(\eta,\mu)\in\bigl( \mathcaligr{P}\bigl(S^{(l)}\bigr)^{\mathbb
{N}} \bigr)
^{2}\qquad \Vert\eta-\mu\Vert:={\sup_{n\geq0}}\Vert\eta_{n}-\mu
_{n}\Vert.
\]
We extend a given integral operator $\mu\in\mathcaligr
{P}(S^{(l)})\mapsto\mu
L\in\mathcaligr{P}(S^{(l+1)})$ into a mapping
\[
\eta=(\eta_{n})_{n\geq0} \in\mathcaligr{P}\bigl(S^{(l)}\bigr)^{\mathbb
{N}}\quad\mapsto\quad
\eta L=(\eta_{n}L)_{n\geq0}\in\mathcaligr{P}\bigl(S^{(l+1)}\bigr)^{\mathbb{N}}.
\]
Sometimes, we slightly abuse the notation and we denote by $\nu$
instead of
$(\nu)_{n\geq0}$ the constant distribution flow equal to a given
measure $%
\nu\in\mathcaligr{P}(S^{(l)})$.

\subsection{Time averaged semigroups}
\label{timeav}

We associate with the mappings $\Phi_{l}$ introduced
in (\ref%
{phi}) the mappings
\[
\Phi^{(l)} \dvtx\eta\in\mathcaligr{P}\bigl(S^{(l-1)}\bigr)^{\mathbb{N}}
\quad\mapsto\quad
\Phi
^{(l)}(\eta)= \bigl( \Phi_{n}^{(l)}(\eta) \bigr) _{n\geq0}\in
\mathcaligr{P}%
\bigl(S^{(l)}\bigr)^{\mathbb{N}}
\]
defined by the coordinate mappings
\[
\forall\eta\in\mathcaligr{P}\bigl(S^{(l-1)}\bigr)^{\mathbb{N}} ,\forall
n\geq0\qquad
\Phi_{n}^{(l)}(\eta):=\Phi_{l}(\eta_{n}).
\]
We denote by
\[
\Phi^{(k,l)}=\Phi^{(k)}\circ\Phi^{(k-1,l)}
\]
with $0\leq l\leq k$, the semigroup associated with the mappings $\Phi
^{(l)}$%
. We also consider the time averaged transformations
\[
\overline{\Phi}{}^{(l)} \dvtx\eta\in\mathcaligr{P}\bigl(S^{(l-1)}\bigr)^{\mathbb
{N}} \quad\mapsto\quad
\overline{\Phi}{}^{(l)}(\eta)= \bigl( \overline{\Phi}_{n}{}^{(l)}(\eta
) \bigr)
_{n\geq0}\in\mathcaligr{P}\bigl(S^{(l)}\bigr)^{\mathbb{N}}
\]
defined by the coordinate mappings
\begin{eqnarray*}
\forall\eta\in\mathcaligr{P}\bigl(S^{(l-1)}\bigr)^{\mathbb{N}},
\forall n\geq0 \qquad
\overline{\Phi}_{n}{}^{(l)}(\eta):\!&=&\frac{1}{n+1}\sum_{p=0}^{n}\Phi
_{p}^{(l)}(\eta) \\
&=&\frac{1}{n+1}\sum_{p=0}^{n}\Phi_{l}(\eta_{p})\in\mathcaligr{P}\bigl(S^{(l)}\bigr).
\end{eqnarray*}
For $l=0$, we use the convention $\Phi_{0}(\eta_{p})=\pi^{(0)}$ for
any $%
0\leq p\leq n$, so that with some abusive but obvious notation
$\overline{%
\Phi}^{(0)}(\eta)=\pi^{(0)}$ represents the constant sequence $%
(\pi^{(0)})_{n\geq0}$ such that $\pi_{n}^{(0)}=\pi^{(0)}$.

We also denote $\overline{\Phi}{}^{(k,l)} \dvtx \mathcaligr{P
}(S^{(l-1)})^{\mathbb{N%
}}\rightarrow\mathcaligr{P }(S^{(k)})^{\mathbb{N}}$ with $0\leq l\leq
k$, the
semigroup associated with the mappings $\overline{\Phi}{}^{(l)}$ and defined
by
\[
\overline{\Phi}{}^{(k,l)}:=\overline{\Phi}{}^{(k)}\circ\overline{\Phi
}%
^{(k-1)}\circ\cdots\circ\overline{\Phi}{}^{(l)}.
\]
We use the convention $\overline{\Phi}{}^{(k,l)}=\mathrm{Id}$, the identity operator,
for $l>k$.

\subsection{Integral operators}

We associate with the kernel $\Gamma_{k}$ from $\mathcaligr{B}(S^{(k)})$
into $%
\mathcaligr{B}(S^{(k-1)})$ introduced in (\ref{lipcn}) the kernel
$\overline{%
\Gamma}^{(k)}$ from $(\mathbb{N}\times\mathcaligr{B}(S^{(k)}))$ into
the set $%
(\mathbb{N}\times\mathcaligr{B}(S^{(k-1)}))$ defined by
%
\begin{eqnarray} \label{defsigma}
\overline{\Gamma}{}^{(k)}((n,f),d(p,g)):=\Sigma(n,dp) \times\Gamma
_{k}(f,dg)\nonumber\\[-8pt]\\[-8pt]
\eqntext{\mbox{with } \Sigma(n,dp):=\dfrac{1}{n+1}%
\displaystyle\sum_{q=0}^{n} \delta_{q}(dp).}
\end{eqnarray}
The semigroup $\overline{\Gamma}{}^{(l_{2},l_{1})}$ ($0\leq l_{1}\leq l_{2}$)
associated with the integral operators $\overline{\Gamma}{}^{(l)}$ is defined
by
\[
\overline{\Gamma}{}^{(l_{2},l_{1})}:=\overline{\Gamma}{}^{(l_{2})}\overline{%
\Gamma}^{(l_{2}-1)}\cdots\overline{\Gamma}{}^{(l_{1})}.
\]
For $l_{1}=l_{2}=0$, we use the convention $\overline{\Gamma}{}^{(0,0)}=
\overline{\Gamma}{}^{(0)}=0$ for the null measure on $(\mathbb
{N}\times
\mathcaligr{B}(S^{(0)}))$. Also observe that
\[
\overline{\Gamma}{}^{(l_{2},l_{1})}=\Sigma^{l_{2}-l_{1}+1}\times
\Gamma
_{l_{2},l_{1}},
\]
where the semigroups $\Sigma^{l_{1}}$ and $\Gamma_{l_{2},l_{1}}$,
$0\leq
l_{1}\leq l_{2}$ associated with the pair of integral operators $\Sigma$
and $\Gamma_{l}$ are
\[
\Sigma^{l_{1}}=\Sigma\Sigma^{l_{1}-1}=\Sigma^{l_{1}-1}\Sigma%
\quad\mbox{and}\quad \Gamma_{l_{2},l_{1}}:=\Gamma_{l_{2}}\Gamma
_{l_{2}-1}\cdots\Gamma_{l_{1}}.
\]
We use the convention $\Sigma^{0}=\mathrm{Id}$.

We end this section with a technical lemma relating the regularity
properties (\ref{lipcn}) of the mappings $\Phi_{k}$ to the regularity
properties of the semigroups $\overline{\Phi}{}^{(k,l)}$.
\begin{lem}
\label{phibar} For any $0\leq l_{1}\leq l_{2}$, $n\geq0$, any flow of
measures $\eta,\mu\in\mathcaligr{P }(S^{(l_{1}-1)})^{\mathbb{N}}$ and any
function $f\in\mathcaligr{B}(S^{(l_{2})})$ we have
\begin{eqnarray*}
&& \bigl\vert\bigl[ \overline{\Phi}{}^{(l_{2},l_{1})}_{n}(\eta
)-\overline{\Phi}%
{}^{(l_{2},l_{1})}_{n}(\mu) \bigr] (f) \bigr\vert
\\
&&\qquad\leq\int_{(\mathbb{N}\times\mathcaligr
{B}(S^{(l_{1}-1)}))} %
\vert[ \eta_{p}-\mu_{p} ] (g) \vert\overline
{\Gamma}{}
^{(l_{2},l_{1})}((n,f),d(p,g)).%
\end{eqnarray*}
\end{lem}
\begin{pf}
Notice that we have $\overline{\Gamma}{}^{(l,l)}=\overline{\Gamma}{}^{(l)}$. We
also observe that $\overline{\Gamma}{}^{(l_{2},l_{1})}$ is a kernel
from $(%
\mathbb{N}\times\mathcaligr{B}(S^{(l_{2})}))$ into $(\mathbb{N}\times
\mathcaligr{B}%
_{n}(S^{(l_{1}-1)}))$. We prove the lemma by induction on the parameter
$%
k=l_{2}-l_{1}$. The result is clearly true for $k=0$. Indeed, by (\ref
{lipcn}%
) we find that for any $l\geq0$
\begin{eqnarray*}
\bigl\vert\bigl[ \overline{\Phi}{}^{(l)}_{n}(\eta)-\overline{\Phi}
{}^{(l)}_{n}(\mu) \bigr] (f) \bigr\vert
&\leq&
\frac{1}{n+1}\sum
_{p=0}^{n} %
\vert[ \Phi_{l}(\eta_{p})-\Phi_{l}(\mu_{p}) ]
(f) \vert\\
&\leq&
\frac{1}{n+1}\sum_{p=0}^{n}\int_{\mathcaligr
{B}(S^{(l-1)})} \vert%
[ \eta_{p}-\mu_{p} ] (g) \vert\Gamma(f,dg).
\end{eqnarray*}
Rewritten in terms of $\overline{\Gamma}{}^{(l)}$, we have proved
that
\[
\bigl\vert\bigl[ \overline{\Phi}{}^{(l)}_{n}(\eta)-\overline{\Phi}
{}^{(l)}_{n}(\mu) \bigr] (f) \bigr\vert\leq\int_{(\mathbb{N}\times
\mathcaligr{B}%
(S^{(l-1)}))} \vert[ \eta_{p}-\mu_{p} ] (g)
\vert%
\overline{\Gamma}{}^{(l)}((n,f),d(p,g)).
\]
This ends the proof of the result for $k=0$. Now, suppose we have proved
that
\[
\bigl\vert\bigl[ \overline{\Phi}{}^{(l_{2},l_{1})}_{p}(\eta
)-\overline{\Phi}{}
^{(l_{2},l_{1})}_{p}(\mu) \bigr] (g) \bigr\vert\leq\int
\vert[
\eta_{q}-\mu_{q} ] (h) \vert\overline{\Gamma}%
{}^{(l_{2},l_{1})}((p,g),d(q,h))
\]
for any pair of integers $l_{1}<l_{2}$ with $l_{2}-l_{1}=k$ for some
$k\geq1$%
. In this case, for any $l<k$ and any function $f\in\mathcaligr{B}(S^{(l+1)})$,
we have
\begin{eqnarray*}
&& \bigl\vert\bigl[ \overline{\Phi}{}^{(l+1,l-k)}_{n}(\eta)-\overline
{\Phi}{}
^{(l+1,l-k)}_{n}(\mu) \bigr] (f) \bigr\vert\\
&&\qquad
= \bigl\vert\bigl[ \overline{\Phi}{}^{(l+1)}_{n}\bigl(\overline{\Phi}{}
^{(l,l-k)}(\eta)\bigr)-\overline{\Phi}{}^{(l+1)}_{n}\bigl(\overline{\Phi}{}
^{(l,l-k)}(\mu)\bigr) \bigr] (f) \bigr\vert%
\end{eqnarray*}
and therefore
\begin{eqnarray*}
&& \bigl\vert\bigl[ \overline{\Phi}{}^{(l+1,l-k)}_{n}(\eta)-\overline
{\Phi}{}
^{(l+1,l-k)}_{n}(\mu) \bigr] (f) \bigr\vert\\
&&\qquad
\leq\int \bigl\vert\bigl[ \overline{\Phi
}{}^{(l,l-k)}_{p}(\eta)-%
\overline{\Phi}{}^{(l,l-k)}_{p}(\mu) \bigr] (g) \bigr\vert
\overline{\Gamma}{}
^{(l+1)}((n,f),d(p,g)).%
\end{eqnarray*}
Under our induction hypothesis, this implies that
\begin{eqnarray*}
&& \bigl\vert\bigl[ \overline{\Phi}{}^{(l+1,l-k)}_{n}(\eta)-\overline
{\Phi}%
{}^{(l+1,l-k)}_{n}(\mu) \bigr] (f) \bigr\vert\\
&&\qquad\leq\int \vert[ \eta_{q}-\mu_{q} ]
(h) \vert \int\overline{\Gamma}{}^{(l+1)}((n,f),d(p,g)) %
\overline{\Gamma}{}^{(l,l-k)}((p,g),d(q,h))
\\
&&\qquad= \int \vert[ \eta_{q}-\mu_{q} ]
(h) \vert%
\overline{\Gamma}{}^{(l+1,l-k)}((n,f),d(q,h)).%
\end{eqnarray*}
Letting $l_{1}=(l-k)$ and $l_{2}=(l+1)$, we have proved that for any $%
l_{1}<l_{2}$ with $l_{2}-l_{1}=(k+1)$
\[
\bigl\vert\bigl[ \overline{\Phi}{}^{(l_{2},l_{1})}_{n}(\eta
)-\overline{\Phi}%
{}^{(l_{2},l_{1})}_{n}(\mu) \bigr] (f) \bigr\vert\leq\int
\vert[
\eta_{p}-\mu_{p} ] (g) \vert\overline{\Gamma}%
{}^{(l_{2},l_{1})}((n,f),d(p,g)).
\]
This ends the proof of the lemma.
\end{pf}

\subsection{Path space semigroups}

To simplify the presentation, we fix a time horizon $m\geq1$ and write
$%
\omega$ instead of $\omega_{K_{\eta}^{(m)}}$, the invariant measure mapping
defined in~(\ref{invpath}). We also write $E$ instead of $E_{m}$.

We extend the mapping $\omega$ on $\mathcaligr{P }(E)$ to $\mathcaligr{P
}(E)^{%
\mathbb{N}}$ by setting
\[
\omega\dvtx\eta=(\eta_{n})_{n\geq0}\in\mathcaligr{P }(E)^{\mathbb
{N}} \quad\mapsto\quad
\omega(\eta)=(\omega_{n}(\eta))_{n\geq0}\in\mathcaligr{P
}(E)^{\mathbb{N}}
\]
with the coordinate mappings $\omega_{n}$ defined by
\[
\omega_{n}(\eta):=\omega(\eta_{n})=\pi^{(0)}\otimes\Phi_{1}\bigl(\eta
^{(0)}_{n}\bigr)%
\otimes\cdots\otimes\Phi_{m}\bigl(\eta^{(m-1)}_{n}\bigr).
\]
For every $l\leq m$, we recall that $\eta^{(l)}_{n}$ stands for the image
measure on $S^{(l)}$ of a given measure $\eta_{n}\in\mathcaligr{P
}(E_{m})$. We
also consider the mappings
\[
\overline{\omega} \dvtx\eta\in\mathcaligr{P }(E)^{\mathbb{N}}
\quad\mapsto\quad
\overline{%
\omega}(\eta)= ( \overline{\omega}_{n}(\eta) ) _{n\geq
0}\in\mathcaligr{%
P }(E)^{\mathbb{N}}
\]
defined by the coordinate mappings
\begin{eqnarray*}
&&\forall\eta=(\eta_{n})_{n\geq0}\in\mathcaligr{P }(E)^{\mathbb
{N}}, \forall n\geq0\\
&&\qquad \overline{\omega}_{n}(\eta):=\frac{1}{n+1}\sum
_{p=0}^{n}\omega
_{p}(\eta)=\frac{1}{n+1}\sum_{p=0}^{n}\omega(\eta_{p}).
\end{eqnarray*}
\begin{lem}
For any $1\leq k\leq m$ and any flow of measures $\eta\in\mathcaligr{P
}(E)^{%
\mathbb{N}}$, we have
\[
\omega^{k}(\eta)=\overline{\pi}{}^{[k-1]}\otimes\bigotimes
_{i=0}^{m-k}%
\Phi^{(i+k,i+1)}\bigl(\eta^{(i)}\bigr).
\]
For $k=m+1$, we have
\[
\forall\eta\in\mathcaligr{P }(E)^{\mathbb{N}}\qquad \omega
^{m+1}(\eta)=\pi^{[m]}.
\]
\end{lem}
\begin{pf}
We use a simple induction on the parameter $k$. The result is clearly true
for $k=1$. Suppose we have proved the result at some rank $k$. In this case
we have
\begin{eqnarray*}
\omega^{k}(\omega(\eta))
&=&
\overline{\pi}{}^{[k-1]}\otimes\Phi_{k,1}\bigl(\omega(%
\eta)^{(0)}\bigr)\otimes\bigotimes_{i=1}^{m-k}\Phi_{i+k,i+1}\bigl(\omega
(\eta)^{(i)}\bigr)
\\
&=&\overline{\pi}{}^{[k-1]}\otimes\pi^{(k)}\otimes\bigotimes
_{i=1}^{m-k}%
\Phi_{i+k,i}\bigl(\eta^{(i-1)}\bigr) \\
&=&\overline{\pi}{}^{[k]}\otimes\bigotimes
_{i=0}^{m-(k+1)}\Phi_{i+(k+1),i+1}\bigl(\eta^{(i)}\bigr).
\end{eqnarray*}
This ends the proof of the lemma.
\end{pf}
\begin{lem}
\label{propp} For any $1\leq k\leq m$ and any $\eta=(\eta
_{n})_{n\geq0}\in%
\mathcaligr{P }(E)^{\mathbb{N}}$, we have
\[
\overline{\omega}{}^{k}_{n}(\eta)=\frac{1}{n+1}\sum_{p=0}^{n} \Biggl[
\overline{%
\pi}^{[k-1]}\otimes\bigotimes_{i=0}^{m-k}\Phi^{(i+k)}_{p} \bigl(
\overline{%
\Phi}^{(i+(k-1),i+1)}\bigl(\eta^{(i)}\bigr) \bigr) \Biggr].
\]
For $k=m+1$, we have
\[
\forall\eta\in\mathcaligr{P }(E)^{\mathbb{N}}\qquad \overline{\omega
}%
^{m+1}(\eta)=\pi^{[m]}.
\]
\end{lem}
\begin{pf}
We use a simple induction on the parameter $k$. The result is clearly true
for $k=1$. Indeed, we have in this case
\[
\overline{\omega}_{n}(\eta)=\frac{1}{n+1}\sum_{p=0}^{n} \Biggl[
\overline{\pi}%
^{[k-1]}\otimes\bigotimes_{i=0}^{m-1}\Phi_{p}^{(i+1)} \bigl( \eta
^{(i)} \bigr) \Biggr].
\]
We also observe that
\[
\overline{\omega}_{n}(\eta)^{(i)}=\frac{1}{n+1}\sum_{p=0}^{n}\Phi
_{p}^{(i)}%
\bigl( \eta^{(i-1)} \bigr) =\overline{\Phi}_{n}{}^{(i)} \bigl(
\eta^{(i-1)} \bigr) \quad\Rightarrow\quad\overline{\omega}(\eta
)^{(i)}=\overline{\Phi}%
^{(i)} \bigl( \eta^{(i-1)} \bigr).
\]

Suppose we have proved the result at some rank $k$. In this case, we have
\[
\overline{\omega}{}^{k}(\overline{\omega}(\eta)) = \frac
{1}{n+1}\sum
_{p=0}^{n} \Biggl[ \overline{\pi}{}^{[k]}\otimes\bigotimes_{i=1}^{m-k}%
\Phi^{(i+k)}_{p} \bigl( \overline{\Phi}{}^{(i+(k-1),i)}\bigl(\eta^{(i-1)}\bigr)
\bigr) %
\Biggr]
\]
from which we conclude that
\[
\overline{\omega}{}^{k+1}(\eta)= \frac{1}{n+1}\sum_{p=0}^{n} \Biggl[
\overline{%
\pi}^{[k]}\otimes\bigotimes_{i=0}^{m-(k+1)}\Phi
{}^{(i+(k+1))}_{p} \bigl(
\overline{\Phi}{}^{(i+k,i+1)}\bigl(\eta^{(i)}\bigr) \bigr) \Biggr].
\]
This ends the proof of the lemma.
\end{pf}

\section{Asymptotic analysis}
\label{asymp}

\subsection{Introduction}

This section is concerned with the asymptotic behavior of i-MCMC models as
the time index $n$ tends to infinity.

The strong law of large numbers is discussed in Section \ref{slln}. We
present nonasymptotic $\mathbb{L}_{r}$-inequalities\vspace*{1pt} that allow us to
quantify the convergence of the occupation measures $\eta
_{n}^{(k)}=\frac{1}{%
n+1}\sum_{p=0}^{n}\delta_{X_{p}^{(k)}}$ of i-MCMC models toward the
solution $\pi^{(k)}$ of the measure-valued equation (\ref{phi}).

Section \ref{unifth} is concerned with uniform convergence results with
respect to the level index $k$. We examine this important question in terms
of the stability properties of the time averaged semigroups introduced in
Section \ref{timeav}. We present nonasymptotic $\mathbb{L}_{r}$%
-inequalities for a series of i-MCMC models that do not depend on the number
of levels. These estimates are probably the most important in practice since
they allow us to quantify the running time of a i-MCMC to achieve a given
precision independently of the time horizon of the limiting measure-valued
equation (\ref{phi}).

Our approach is based on an original combination of nonlinear semigroup
techniques with the asymptotic analysis of time inhomogeneous Markov chains
developed in Section \ref{timc}. The following technical lemma
presents a
more or less well-known generalized Minkowski integral inequality which will
be used in our proofs. 
\begin{lem}[(Generalized Minkowski integral inequality)]
\label{mink} For any pair of bounded positive measures $\mu_{1}$ and $
\mu_{2} $ on some measurable spaces $(E_{1}, \mathcaligr{E }_{1})$ and $(E_{2},
\mathcaligr{E }_{2})$, any bounded measurable function $\varphi$ on the product
space $(E_{1}\times E_{2})$ any $p\geq1$, we have
\begin{eqnarray*}
&&
\biggl[ \int_{E_{1}} \mu_{1}(dx_{1}) \biggl|
\int_{E_{2}} \varphi(x_{1},x_{2}) \mu_{2}(dx_{2}) \biggr|
^{p} \biggr] ^{%
{1/p}}
\\
&&\qquad\leq\int_{E_{2}} \biggl(
\int_{E_{1}} |\varphi(x_{1},x_{2})|^{p} \mu_{1}(dx_{1}) \biggr)
^{1/p} \mu_{2}(dx_{2}).%
\end{eqnarray*}
\end{lem}
\begin{pf}
Without loss of generality, we suppose that $\varphi$ is a nonnegative
function. For $p=1$, the lemma is a direct consequence of Fubini's theorem.
Let us assume that $p>1$, and let $p^{\prime}$ be such that $\frac{1}{
p^{\prime}}+\frac{1}{p}=1$. First, we notice that the functions
\[
\varphi_{1}(x_{1}):=\int_{E_{2}} \varphi(x_{1},x_{2}) \mu
_{2}(dx_{2}) %
\quad\mbox{and}\quad \phi_{p}(x_{2}):= \biggl( \int_{E_{1}} |\varphi
(x_{1},x_{2})|^{p} \mu_{1}(dx_{1}) \biggr) ^{1/p}
\]
are measurable for every $p\geq1$. In this notation, we need to prove
that $%
\mu_{1}(\varphi_{1}^{p})^{1/p}\leq\mu_{2}(\phi_{p})$. It
is also
convenient to consider the function
\[
\psi(x_{1},x_{2}):=\varphi(x_{1},x_{2})/\phi_{p}(x_{2})^{
{1}/{p^{\prime}}.%
}
\]
We use the convention $\psi(x_{1},x_{2})=0$, for every $x_{1}\in
E_{1}$ as
long as $\phi_{p}(x_{2})=0$. We observe that
\[
\biggl( \int_{E_{1}} \psi(x_{1},x_{2})^{p} \mu_{1}(dx_{1}) \biggr)
^{1/p}={\phi_{p}(x_{2})}/{\phi_{p}(x_{2})^{{1}/{p^{\prime}}}}=\phi
_{p}(x_{2})^{{1}/{p}}.
\]
By construction, we have
\begin{eqnarray*}
\varphi_{1}(x_{1}) & = & \int_{E_{2}} \psi(x_{1},x_{2}) \phi
_{p}(x_{2})^{{1%
}/{p^{\prime}}} \mu_{2}(dx_{2}) \\
& \leq & \biggl[ \int_{E_{2}} \psi(x_{1},x_{2})^{p} \mu
_{2}(dx_{2}) \biggr] ^{%
{1/p}}\times\mu_{2}(\phi_{p})^{{1}/{p^{\prime}}}
\end{eqnarray*}
from which we conclude that
\begin{eqnarray*}
\mu_{1}(\varphi_{1}^{p})
&\leq&\mu_{2}(\phi_{p})^{
{p}/{p^{\prime}}}\times%
\biggl[ \int_{E_{2}} \psi(x_{1},x_{2})^{p} \mu_{1}(dx_{1})\mu
_{2}(dx_{2})%
\biggr] \\
&=&\mu_{2}(\phi_{p})^{{p}/{p^{\prime}}}\times\mu_{2}(\phi
_{p})=\mu_{2}(%
\phi_{p})^{p}.
\end{eqnarray*}
The end of the proof is now clear.
\end{pf}

\subsection{Strong law of large numbers}
\label{slln}

This section is mainly concerned with the proof of the following
$\mathbb{L}%
_{r}$-inequalities for the occupation measure of an i-MCMC model at a given
level.
\begin{theorem}
Under the regularity\vspace*{1pt} conditions (\ref{lipcn}) and (\ref{condin}), we have
for any $k\geq0$, any function $f\in\mathcaligr{B}_{1}(S^{(k)})$ and
any $%
n\geq0 $ and $r\geq1$
%
\begin{eqnarray}\label{eqf}
&&\sqrt{(n+1)} \mathbb{E} \bigl( \bigl\vert\bigl[ \eta
_{n}^{(k)}-\pi^{(k)}%
\bigr] (f) \bigr\vert^{r} \bigr)^{1/r}
\nonumber\\[-8pt]\\[-8pt]
&&\qquad
\leq e(r) \sum
_{l=0}^{k}(1+c_{l}) \biggl( \frac{n_{l}}{1-b_{l}(n_{l})} \biggr)
^{2}\prod_{l+1\leq i\leq k}2\Lambda_{i}.\nonumber
\end{eqnarray}
\end{theorem}
\begin{pf}
We prove the theorem by induction on the parameter $k$. First, we observe
that the estimate (\ref{eqf}) is true for $k=0$. Indeed, by
Corollary \ref%
{corinter} we have that
\[
\sqrt{(n+1)} \mathbb{E} \bigl( \bigl\vert\bigl[ \eta
_{n}^{(0)}-\pi^{(0)}%
\bigr] (f) \bigr\vert^{r} \bigr) ^{{1}/{r}}\leq
e(r) (1+c_{0}) \biggl(
\frac{n_{0}}{1-b_{0}(n_{0})} \biggr) ^{2}%
\]
for some finite constant $e(r)<\infty$ whose value only depends on the
parameter $r$. We further suppose that the estimate (\ref{eqf}) is
true at
rank $(k-1)$. To prove that it is also true at rank $k$, we use the
decomposition
%
\begin{equation} \label{induc}\qquad
\bigl[ \eta{}^{(k)}_{n}-\pi^{(k)} \bigr] = \bigl[ \eta
{}^{(k)}_{n}-\overline{\Phi}%
{}^{(k)}_{n}\bigl(\eta^{(k-1)}\bigr) \bigr] + \bigl[ \overline{\Phi}{}^{(k)}_{n}\bigl(%
\eta{}^{(k-1)}\bigr)-\overline{\Phi}{}^{(k)}_{n}\bigl(\pi{}^{(k-1)}\bigr) \bigr].
\end{equation}
For every $k\geq0$, given a realization of the chain $%
X^{(k-1)}:=(X_{p}^{(k-1)})_{p\geq0}$ the $k$th level chain $X_{n}^{(k)}$
behaves as a Markov chain with random Markov transitions $%
M_{\eta_{n}^{(k-1)}}^{(k)}$ dependent on the current occupation
measure of
the chain at level $(k-1)$. Therefore, using Corollary \ref{corinter} again
we notice that
\[
\sqrt{(n+1)} \mathbb{E} \bigl( \bigl\vert\bigl[ \eta
{}^{(k)}_{n}-\overline{\Phi%
}{}^{(k)}_{n}\bigl(\eta^{(k-1)}\bigr) \bigr] (f) \bigr\vert^{r} \bigr) ^{
{1}/{r}%
}\leq e(r) (1+c_{k}) \biggl( \frac{n_{k}}{1-b_{k}(n_{k})} \biggr) ^{2}%
\]
for some finite constant $e(r)<\infty$ whose values only depends on the
parameter $r$.

Using the decomposition (\ref{induc}) and Lemma \ref{phibar}, we obtain
\begin{eqnarray*}
&& \bigl\vert\bigl[ \eta_{n}^{(k)}-\pi^{(k)} \bigr]
(f) \bigr\vert\\
&&
\qquad\leq\bigl\vert\bigl[ \eta_{n}^{(k)}-\overline{\Phi
}{}^{(k)}_{n}\bigl(\eta^{(k-1)}\bigr)%
\bigr] (f) \bigr\vert
\\
&&\qquad\quad{} + \int \bigl\vert\bigl[
\eta_{p}^{(k-1)}-\pi^{(k-1)} \bigr] (g) \bigr\vert\overline
{\Gamma}%
{}^{(k)}((n,f),d(p,g)).%
\end{eqnarray*}
For every function $f\in\mathcaligr{B}_{1}(S^{(l)})$, and any $n\geq0$,
$k\geq0$%
, $r\geq1$, we set
\[
J_{n}^{(k)}(f):=\sqrt{n+1} \mathbb{E} \bigl( \bigl\vert\bigl[
\eta_{n}^{(k)}-\pi^{(k)} \bigr] (f) \bigr\vert^{r} \bigr) ^{
{1}/{r}%
} \quad\mbox{and}\quad j^{(k)}:=\sup_{n\geq1} \sup_{f\dvtx\Vert
f\Vert\leq
1}J_{n}^{(k)}(f).
\]
By the generalized Minkowski integral inequality presented in
Lemma \ref%
{mink}, we find that
\begin{eqnarray*}
J_{n}^{(k)}(f)&\leq& e(r)(1+c_{k}) \biggl( \frac
{n_{k}}{1-b_{l}(n_{k})}%
\biggr) ^{2} \\
&&{} +\sqrt{n+1}\int J_{p}^{(k-1)}(g) \frac{1}{\sqrt
{p+1}} \overline{%
\Gamma}{}^{(k)}((n,f),d(p,g)).
\end{eqnarray*}
Since we have
%
\begin{equation} \label{sqrt2}
\int_{\mathbb{N}} \frac{1}{\sqrt{q+1}} \Sigma
(n,dq)=\frac{1}{n+1%
}\sum_{q=0}^{n}\frac{1}{\sqrt{q+1}}\leq\frac{2}{\sqrt{n+1}}
\end{equation}
we conclude that
\[
J_{n}^{(k)}(f)\leq e(r)(1+c_{k}) \biggl( \frac
{n_{k}}{1-b_{l}(n_{k})} \biggr)
^{2}+2 j^{(k-1)} \sup_{f}\int\Vert g\Vert\Gamma_{k}(f,dg)
\]
and therefore
\[
j^{(k)}\leq e(r)(1+c_{k}) \biggl( \frac{n_{k}}{1-b_{k}(n_{k})} \biggr)
^{2}+j^{(k-1)} 2\Lambda_{k} .
\]
Under the induction hypothesis, we have
\[
j^{(k-1)} 2\Lambda_{k} \leq e(r) \sum_{l=0}^{k-1}(1+c_{l}) \biggl(
\frac{%
n_{l}}{1-b_{l}(n_{l})} \biggr) ^{2}\prod_{l+1\leq i\leq k}2\Lambda_{i}
\]
and therefore
\begin{eqnarray*}
j^{(k)} & \leq& e(r) \biggl[ (1+c_{k}) \biggl( \frac
{n_{k}}{1-b_{k}(n_{k})}%
\biggr) ^{2}\\
&&\hspace*{22.22pt}{} + \sum_{l=0}^{k-1}(1+c_{l}) \biggl( \frac
{n_{l}}{%
1-b_{l}(n_{l})} \biggr) ^{2}\prod_{l+1\leq i\leq k}2\Lambda_{i}
\biggr] \\
& =&\sum_{l=0}^{k}(1+c_{l}) \biggl( \frac
{n_{l}}{1-b_{l}(n_{l})} \biggr)
^{2}\prod_{l+1\leq i\leq k}2\Lambda_{i}.
\end{eqnarray*}
This ends the proof of the theorem.
\end{pf}

\subsection{A uniform convergence theorem}
\label{unifth}

This section focuses on the behavior of an i-MCMC model associated with a
large number of levels. We establish an uniform convergence theorem under
the assumption that the time averaged semigroup $\overline{\Phi}{}^{(k,l)}$
introduced in Section \ref{timeav} is exponentially stable; that is, there
exist some positive constants $\lambda_{1},\lambda_{2}>0$ and an
integer $%
k_{0}$ such that for every $l\geq0$, $\eta,\mu\in\mathcaligr
{P}(S^{(l)})^{%
\mathbb{N}}$ and any $k\geq k_{0}$ we have
%
\begin{equation} \label{exponentialstability}
\bigl\Vert\overline{\Phi}{}^{(l+k,l+1)}(\eta)-\overline{\Phi}{}^{(l+k,l+1)}(\mu
)\bigr\Vert\leq\lambda_{1} e^{-\lambda_{2}k}.
\end{equation}
We also assume that the parameters $(b_{k},c_{k},n_{k},\Lambda_{k})$ are
chosen so that
%
\begin{equation}\label{exponentialstability2}\quad
A=\sup_{k\geq0} \biggl[ (1+c_{k}) \biggl( \frac
{n_{k}}{1-b_{k}(n_{k})} \biggr)
^{2} \biggr] <\infty\quad\mbox{and}\quad B:={2}\sup_{k\geq
1}{\Lambda
_{k}}<\infty.
\end{equation}
For the Feynman--Kac transformations (\ref{fksg}), we give in Section
\ref%
{FKS} sufficient conditions on ${G_{l}}$ and ${L_{l+1}}$ ensuring (\ref
{exponentialstability}) is satisfied. If (\ref{exponentialstability})
and (%
\ref{exponentialstability2}) are both satisfied, we have the following
uniform convergence result:
\begin{theorem}
If $B=1$, then we have for any $r\geq1$, any parameter $n$ such that $%
(n+1)\geq e^{2\lambda_{2}(k_{0}+1)}$, and for any $(f_{l})_{l\geq
0}\in%
\prod_{l\geq0}\operatorname{Osc}_{1}(S^{(l)})$
\[
\sup_{l\geq0}{\mathbb{E} \bigl( \bigl| \bigl[ \eta_{n}^{(l)}-\pi
^{(l)} \bigr]
(f_{l}) \bigr| ^{r} \bigr) ^{{1}/{r}}}\leq\frac{e(r)}{\sqrt
{n+1}} \biggl(
A \biggl( 1+\frac{\log{(n+1)}}{2\lambda_{2}} \biggr)
+\lambda_{1} e^{\lambda_{2}} \biggr).
\]
If $B>1$, then we have for any $r\geq1$, any $n$ such that $(n+1)\geq
e^{2(\lambda_{2}+\log{B})(k_{0}+1)}$, and for any $(f_{l})_{l\geq
0}\in
\prod_{l\geq0}\operatorname{Osc}_{1}(S^{(l)})$.
\[
\sup_{l\geq0}{\mathbb{E} \bigl( \bigl| \bigl[ \eta_{n}^{(l)}-\pi
{}^{(l)} \bigr]
(f_{l}) \bigr| ^{r} \bigr) ^{{1}/{r}}}\leq e(r) \biggl[ \frac
{AB}{B-1}%
+\lambda_{1} \biggr] \frac{e^{\lambda_{2}}}{(n+1)^{\alpha/2}}
\]
with $\alpha:=\frac{\lambda_{2}}{(\lambda_{2}+\log{B})}$.
\end{theorem}
\begin{pf}
First, we notice that we have the following estimate from (\ref{eqf})
and (%
\ref{exponentialstability2}) for any $k\geq0$:
%
\begin{equation} \label{eqq1}
\sqrt{(n+1)} \mathbb{E} \bigl( \bigl\vert\bigl[ \eta
_{n}^{(k)}-\pi^{(k)}%
\bigr] (f_{k}) \bigr\vert^{r} \bigr) ^{{1}/{r}}\leq
e(r) A \frac{%
B^{k+1}-1}{B-1}.
\end{equation}
For $B=1$, we use the convention $\frac{B^{k}-1}{B-1}=k$.

We have the following decomposition:
%
\begin{eqnarray} \label{eqq2}\qquad
\eta_{n}^{(l+k)}-\pi^{(l+k)} & = & \bigl[ \eta_{n}^{(l+k)}-\overline
{\Phi}%
{}^{(l+k,l+1)_{n}}\bigl(\eta^{(l)}\bigr) \bigr]\nonumber\\
&&{} + \bigl[ \overline{\Phi}%
{}^{(l+k,l+1)}_{n}\bigl(\eta^{(l)}\bigr)-\overline{\Phi}_{n}{}^{(l+k,l+1)}\bigl(\pi
^{(l)}\bigr)%
\bigr] \nonumber\\[-8pt]\\[-8pt]
& = &\sum_{i=l+1}^{l+k} \bigl[ \overline{\Phi}{}^{(l+k,i+1)}_{n}\bigl(\eta
^{(i)}\bigr)-%
\overline{\Phi}_{n}{}^{(l+k,i+1)}\bigl(\overline{\Phi}{}^{(i)}\bigl(\eta
^{(i-1)}\bigr)\bigr) \bigr]
\nonumber\\
&&{} + \bigl[ \overline{\Phi}%
{}^{(l+k,l+1)}_{n}\bigl(\eta^{(l)}\bigr)-\overline{\Phi}{}^{(l+k,l+1)}_{n}\bigl(\pi
^{(l)}\bigr) \bigr].\nonumber
\end{eqnarray}
Recall that we use the convention $\overline{\Phi}{}^{(l_{1},l_{2})}=\mathrm{Id}$ for $%
l_{1}<l_{2}$, so that
\[
i=l+k\quad\Longrightarrow\quad\overline{\Phi}{}^{(l+k,i+1)}_{n}\bigl(\eta
^{(i)}\bigr)=\overline{%
\Phi}{}^{(l+k,l+k+2)}_{n}\bigl(\eta^{(l+k)}\bigr)=\eta_{n}^{(l+k)}.
\]
Using Lemma \ref{phibar}, we find that
\begin{eqnarray*}
&&\bigl\vert\bigl[ \overline{\Phi}{}^{(l_{2},l_{1}+1)}_{n}\bigl(\eta
^{(l_{1})}\bigr)-%
\overline{\Phi}_{n}{}^{(l_{2},l_{1})}\bigl(\overline{\Phi}{}^{(l_{1})}\bigl(\eta
^{(l_{1}-1)}\bigr)\bigr) \bigr] (f_{l_{2}}) \bigr\vert\\
&&\qquad\leq\int \bigl\vert\bigl[ \eta
_{p}^{(l_{1})}-\overline{\Phi}{}
_{p}^{(l_{1})}\bigl(\eta^{(l_{1}-1)}\bigr) \bigr] (g) \bigr\vert\overline
{\Gamma}%
{}^{(l_{2},l_{1}+1)}((n,f_{l_{2}}),d(p,g)).%
\end{eqnarray*}
By the generalized Minkowski integral inequality, this implies that
\begin{eqnarray*}
&&\mathbb{E} \bigl( \bigl\vert\bigl[ \overline{\Phi
}{}^{(l_{2},l_{1}+1)}_{n}\bigl(%
\eta^{(l_{1})}\bigr)-\overline{\Phi}{}^{(l_{2},l_{1}+1)}_{n}\bigl(\overline
{\Phi}%
{}^{(l_{1})}\bigl(\eta^{(l_{1}-1)}\bigr)\bigr) \bigr] (f_{l_{2}}) \bigr\vert
^{r} \bigr) ^{%
{1}/{r}}
\\
&&\qquad\leq\int\mathbb{E} \bigl(\bigl \vert\bigl[ \eta
{}^{(l_{1})}_{p}-%
\overline{\Phi}{}^{(l_{1})}_{p}\bigl(\eta^{(l_{1}-1)}\bigr) \bigr] (g)
\bigr\vert
^{r} \bigr) ^{{1}/{r}} \overline{\Gamma}%
{}^{(l_{2},l_{1}+1)}((n,f_{l_{2}}),d(p,g)).%
\end{eqnarray*}
Using Corollary \ref{corinter}, we find that
\begin{eqnarray*}
&&\mathbb{E} \bigl( \bigl\vert\bigl[ \overline{\Phi
}{}^{(l_{2},l_{1}+1)}_{n}\bigl(%
\eta^{(l_{1})}\bigr)-\overline{\Phi}{}^{(l_{2},l_{1}+1)}_{n}\bigl(\overline
{\Phi}%
{}^{(l_{1})}\bigl(\eta^{(l_{1}-1)}\bigr)\bigr) \bigr] (f_{l_{2}}) \bigr\vert
^{r} \bigr) ^{%
{1}/{r}}
\\
&&\qquad \leq e(r) (1+c_{l_{1}}) \biggl( \frac
{n_{l_{1}}}{1-b_{l_{1}}(n_{l_{1}})}%
\biggr) ^{2}
\\
&&\qquad\quad{}\times\int_{\{0,\ldots,n\}} \frac{1}{\sqrt{(p+1)}}%
\Sigma^{(l_{2}-l_{1})}(n,dp)\times\int\Vert g\Vert
\Gamma_{l_{2},l_{1}+1}(f_{l_{2}},dg).%
\end{eqnarray*}
By (\ref{sqrt2}) and
\[
\int\Gamma_{k,l}(f_{l_{2}},dg) \Vert g\Vert\leq
\Lambda
_{k,l} \Vert f_{l_{2}}\Vert\qquad\mbox{with } \Lambda_{k,l}\leq\prod_{l\leq i\leq k}\Lambda_{i}\leq
B^{k-l+1}<\infty,
\]
we conclude that
%
\begin{eqnarray}\label{Laref}
&&\sqrt{(n+1)}\mathbb{E} \bigl( \bigl\vert\bigl[ \overline{\Phi}%
{}^{(l_{2},l_{1}+1)}_{n}\bigl(\eta^{(l_{1})}\bigr)-\overline{\Phi
}{}^{(l_{2},l_{1}+1)}_{n}\bigl(%
\overline{\Phi}{}^{(l_{1})}\bigl(\eta^{(l_{1}-1)}\bigr)\bigr) \bigr]
(f_{l_{2}}) \bigr\vert
^{r} \bigr) ^{{1}/{r}} \nonumber\\[-8pt]\\[-8pt]
&&\qquad\leq e(r) A B^{l_{2}-l_{1}} \Vert f_{l_{2}}\Vert.\nonumber
\end{eqnarray}
Using the decomposition (\ref{eqq2}), we prove that for every
$f_{l+k}\in%
\mathcaligr{B}_{1}(S^{(l+k)})$ and any $k\geq k_{0}$
\[
\sup_{l\geq0}{\mathbb{E} \bigl( \bigl\vert\bigl[ \eta
_{n}^{(l+k)}-\pi^{(l+k)}%
\bigr] (f_{l+k}) \bigr\vert^{r} \bigr) ^{{1}/{r}}}\leq
e(r) \frac{A}{%
\sqrt{n+1}} \frac{B^{k}-1}{B-1}+\lambda_{1} e^{-\lambda_{2} k}.
\]
Finally, by (\ref{eqq1}), we conclude that for every $k\geq k_{0}$
\[
\sup_{l\geq0}{\mathbb{E} \bigl( \bigl\vert\bigl[ \eta
_{n}^{(l)}-\pi^{(l)}%
\bigr] (f_{l}) \bigr\vert^{r} \bigr) ^{{1}/{r}}}\leq
e(r) \frac{A}{%
\sqrt{n+1}} \frac{B^{k+1}-1}{B-1}+\lambda_{1} e^{-\lambda_{2} k}.
\]
For $B=1$, we have
\[
\sup_{l\geq0}{\mathbb{E} \bigl( \bigl\vert\bigl[ \eta
_{n}^{(l)}-\pi^{(l)}%
\bigr] (f_{l}) \bigr\vert^{r} \bigr) ^{{1}/{r}}}\leq
e(r) A \frac{(k+1)%
}{\sqrt{n+1}}+\lambda_{1} e^{-\lambda_{2} k}.
\]
In this situation, we choose the parameters $k,n$ such that
\[
k=k(n):=\biggl\lfloor\frac{\log{(n+1)}}{2\lambda_{2}}\biggr\rfloor\geq k_{0}.
\]
Notice that $k(n)$ is the largest integer $k$ satisfying
\[
k\leq\frac{\log{(n+1)}}{2\lambda_{2}}\quad\Leftrightarrow\quad \biggl(
\frac{1}{%
\sqrt{n+1}}\leq e^{-\lambda_{2} k} \biggr).
\]
Since $(k(n)+1)\geq\frac{\log{(n+1)}}{2\lambda_{2}}$, we have
\[
e^{-\lambda_{2} k(n)}\leq e^{\lambda_{2}} e^{-\lambda_{2} (
{\log{(n+1)}})/({%
2\lambda_{2}})}=\frac{e^{\lambda_{2}}}{\sqrt{n+1}}
\]
from which we conclude that
\[
A \frac{(k(n)+1)}{\sqrt{n+1}}+\lambda_{1} e^{-\lambda_{2} k(n)}\leq
\frac{1}{%
\sqrt{n+1}} \biggl( A \biggl( 1+\frac{\log{(n+1)}}{2\lambda
_{2}} \biggr)
+\lambda_{1} e^{\lambda_{2}} \biggr).
\]
For $B>1$, we choose the parameters $k,n$ such that
\[
k=k(n):=\biggl\lfloor\frac{\log{(n+1)}}{2(\lambda_{2}+\log{B})}\biggr\rfloor
\geq k_{0}.
\]
Notice that $k(n)$ is the largest integer $k$ such that
\[
k\leq\frac{\log{(n+1)}}{2(\lambda_{2}+\log{B})}
\quad\Leftrightarrow\quad
\biggl(\frac{B^{k}}{\sqrt{n+1}}\leq e^{-\lambda_{2} k} \biggr).
\]
Since $(k(n)+1)\geq\frac{\log{(n+1)}}{2(\lambda_{2}+\log{B})}$, we have
\[
\frac{B^{k(n)}}{\sqrt{n+1}}\leq e^{-\lambda_{2} k(n)}\leq
e^{\lambda_{2}} e^{-\lambda_{2} ({\log{(n+1)}})/({2(\lambda
_{2}+\log{B})})}=%
\frac{e^{\lambda_{2}}}{(n+1)^{\alpha/2}}
\]
with $\alpha:=\frac{\lambda_{2}}{(\lambda_{2}+\log{B})}$, from
which we
conclude that
\[
\frac{A}{\sqrt{n+1}} \frac{B^{k(n)+1}-1}{B-1}+\lambda
_{1} e^{-\lambda
_{2} k(n)}\leq\biggl[ \frac{AB}{B-1}+\lambda_{1} \biggr] \frac
{e^{\lambda
_{2}}}{(n+1)^{\alpha/2}}-\frac{AB}{B-1} \frac{1}{\sqrt{n+1}}.
\]
This ends the proof of the theorem.
\end{pf}

\section{Path space models}
\label{pathmodels}

In the previous section, we have established $\mathbb{L}_{r}$-mean error
bounds and exponential estimates quantifying the convergence of the
occupation measures $\eta_{n}^{ ( k ) }$ toward the
solutions $\pi
_{n}^{ ( k ) }$ of the measure-valued equation (\ref{phi}).
We show
here that it is also possible to establish such results to quantify the
convergence of the path-space occupation measures $\overline{\eta}_{n}^{[m]}$
introduced in (\ref{pathocc}) toward the tensor product measure
$\overline{%
\pi}^{(m)}$ defined in (\ref{eq:tensorproducttargetmeasure}).

\subsection{$\mathbb{L}_{r}$-mean error bounds}
\label{lrmean}

Our main result is the following theorem:
\begin{theorem}
\label{theo411} For every $f\in\mathcaligr{B }(E_{m})$, we have
\[
\sup_{n\geq1}\sqrt{n}\, \mathbb{E} \bigl( \bigl| \bigl[ \overline
{\eta}%
^{[m]}_{n}-\overline{\pi}{}^{(m)} \bigr] (f) \bigr| ^{r} \bigr)
^{{1}/{r}}%
<\infty.
\]
\end{theorem}
\begin{pf}
To simplify the presentation, we fix a time horizon $m\geq1$ and write
$%
\omega$ instead of $\omega_{K_{\eta}^{(m)}}$, the invariant\vspace*{-2pt} measure
mapping defined in (\ref{invpath}). We also write $E$ instead of $E_{m}$,
and $\overline{\eta}_{n}$ instead of $\overline{\eta}_{n}^{[m]}$.
In this
notation, $(\overline{\eta}^{(l)})$ represents the sequence of occupation
measures $\overline{\eta}_{n}^{(l)}:=\frac{1}{n+1}\sum
_{p=0}^{n}\delta
_{X_{p}^{(l)}}\in\mathcaligr{P}(S^{(l)})$ of the i-MCMC model on the $l$th
level space $S^{(l)}$.

Using the fact that $\overline{\omega}{}^{m+1}(\eta)=\overline{\pi}{}^{[m]}$, we
obtain the following decomposition for any $\eta\in\mathcaligr
{P}(E)^{\mathbb{N}%
}$
%
\begin{equation} \label{decomp1}
\eta-\overline{\pi}{}^{[m]}=\sum_{k=0}^{m} \bigl[ \overline{\omega}{}^{k}(\eta)-%
\overline{\omega}{}^{k+1}(\eta) \bigr].
\end{equation}
In the above-displayed formula, $\overline{\pi}{}^{[m]}=(\overline{\pi
}%
^{[m]}_n)_{n\in\mathbb{N}}\in\mathcaligr{P}(E)^{\mathbb{N}}$ stands
for the
constant sequence of measures $\overline{\pi}{}^{[m]}_n=\overline{\pi}{}^{[m]}$,
for any $n\in\mathbb{N}$.

Using Proposition \ref{propp}, the $k$th iterate $\overline{\omega}{}^{k}$ of
the mapping $\overline{\omega}$ can be rewritten for any $\eta\in
\mathcaligr{P}%
(E)^{\mathbb{N}}$ in the following form:
\[
\overline{\omega}_{n}^{k}(\eta)=\frac{1}{n+1}\sum_{p=0}^{n} \bigl[
\overline{%
\pi}^{[k-1]}\otimes\Pi_{p}^{(k,m)}\bigl(\bigl(\eta^{(l)}\bigr)_{0\leq l\leq
m}\bigr) \bigr].
\]
Here the mappings
\[
\Pi^{(k,m)} \dvtx\mu\in\prod_{0\leq i\leq m}\mathcaligr
{P}\bigl(S^{(i)}\bigr)^{\mathbb{N}%
}\mapsto\Pi^{(k,m)}(\mu)= \bigl( \Pi_{n}^{(k,m)}(\mu) \bigr)
_{n\geq0}\in\Biggl( \bigotimes_{i=k}^{m}\mathcaligr{P}\bigl(S^{(i)}\bigr)
\Biggr) ^{\mathbb{%
N}}
\]
are defined for any $n\geq0$ by
\[
\Pi_{n}^{(k,m)}(\mu):=\bigotimes_{i=0}^{m-k}\Pi_{n}^{(k,m),(i)}(\mu
)\in\bigotimes_{i=0}^{m-k}\mathcaligr{P}\bigl(S^{(i+k)}\bigr)
\]
with for any $(\mu^{(l)})_{0\leq l\leq m}\in\prod_{0\leq i\leq
m}\mathcaligr{P}%
(S^{(i)})^{\mathbb{N}}$ and any $0\leq i\leq m-k$
\[
\Pi_{n}^{(k,m),(i)}\bigl(\bigl(\mu^{(l)}\bigr)_{l}\bigr):=\Phi_{i+k} \bigl( \overline
{\Phi}%
{}^{(i+(k-1),i+1)}_{n}\bigl(\mu^{(i)}\bigr) \bigr) \in\mathcaligr{P}\bigl(S^{(i+k)}\bigr).
\]
We emphasize that $\Pi_{n}^{(k,m)}(\mu)$ only depends on the flow of
measures $(\mu^{(l)})_{0\leq l\leq m-k}$, and
\begin{eqnarray*}
&&\overline{\omega}_{n}^{k+1}(\eta)
\\
&&\qquad= \frac{1}{n+1}\sum
_{p=0}^{n} \bigl[
\overline{\pi}{}^{[k]}\otimes\Pi_{p}^{(k+1,m)}\bigl(\bigl(\eta
^{(l)}\bigr)_{l}\bigr) \bigr] \\
&&\qquad= \frac{1}{n+1}\sum_{p=0}^{n} \Biggl[ \overline{\pi}{}^{[k-1]}\otimes
\pi
^{(k)}\otimes\bigotimes_{i=0}^{m-(k+1)}\Phi_{i+k+1} \bigl( \overline
{\Phi}%
{}^{(i+k,i+2)}_{p}\bigl(\overline{\Phi}{}^{(i+1)}\bigl(\eta^{(i)}\bigr)\bigr) \bigr)
\Biggr] \\
&&\qquad= \frac{1}{n+1}\sum_{p=0}^{n} \Biggl[ \overline{\pi}{}^{[k-1]}\otimes
\bigotimes_{i=0}^{m-k}\Phi_{i+k} \bigl( \overline{\Phi
}{}^{(i+(k-1),i+1)}_{p}\bigl(%
\overline{\Phi}{}^{(i)}\bigl(\eta^{(i-1)}\bigr)\bigr) \bigr) \Biggr]
\end{eqnarray*}
with the convention $\overline{\Phi}{}^{(0)}(\eta^{(-1)}))=\pi
^{(0)}$, for $%
i=0 $. This implies that for any $0\leq k\leq m$
\[
\overline{\omega}_{n}^{k+1}(\eta)=\frac{1}{n+1}\sum_{p=0}^{n} \bigl[
\overline{\pi}{}^{[k-1]}\otimes\Pi_{p}^{(k,m)}\bigl(\bigl(\overline{\Phi}%
{}^{(l)}\bigl(\eta^{(l-1)}\bigr)\bigr)_{l}\bigr) \bigr]
\]
and therefore
%
\begin{eqnarray}\label{yenamarre}\qquad
&&\overline{\omega}_{n}^{k}(\eta)-\overline{\omega}{}^{k+1}_{n}(\eta)
\nonumber\\
&&\qquad=\frac{1}{n+1}\sum_{p=0}^{n} \bigl[ \overline{\pi}{}^{[k-1]}\otimes
\bigl\{
\Pi_{p}^{(k,m)}\bigl(\bigl(\eta^{(l)}\bigr)_{l}\bigr)-\Pi_{p}^{(k,m)} \bigl( \bigl(
\overline{\Phi%
}{}^{(l)}\bigl(\eta^{(l-1)}\bigr) \bigr) _{l} \bigr) \bigr\} \bigr].%
\end{eqnarray}
Moving one step further, we introduce the decomposition
%
\begin{eqnarray}\label{Pi}
&&
\Pi^{(k,m)}(\mu)-\Pi^{(k,m)}(\nu)
\nonumber\\
&&\qquad= \sum_{j=0}^{m-k} \Biggl\{ \Biggl(
\bigotimes_{i=0}^{j-1}\Pi^{(k,m),(i)}(\nu) \Biggr)\nonumber\\[-8pt]\\[-8pt]
&&\qquad\quad\hspace*{22pt}{} \otimes\bigl[
\Pi^{(k,m),(j)}(\mu)-\Pi^{(k,m),(j)}(\nu) \bigr] \nonumber\\
&&\qquad\quad\hspace*{60.4pt}{}\otimes\Biggl(
\bigotimes_{i=j+1}^{m-k}\Pi^{(k,m),(i)}(\mu) \Biggr) \Biggr\}\nonumber
\end{eqnarray}
for any $\mu=(\mu^{(l)})_{0\leq l\leq m}$ and $\nu=(\nu
^{(l)})_{0\leq l\leq
m}\in\prod_{0\leq i\leq m}\mathcaligr{P}(S^{(i)})^{\mathbb{N}}$, with
the flow
of signed measures
\begin{eqnarray*}
&&\Pi_{n}^{(k,m),(j)}(\mu)-\Pi_{n}^{(k,m),(j)}(\nu)
\\
&&\qquad= \bigl[ \Phi_{j+k} \bigl( \overline{\Phi}{}^{(j+(k-1),j+1)}_{n}\bigl(\mu
^{(j)}\bigr)%
\bigr) -\Phi_{j+k} \bigl( \overline{\Phi}{}^{(j+(k-1),j+1)}_{n}\bigl(\nu
^{(j)}\bigr)%
\bigr) \bigr].
\end{eqnarray*}
For every $f\in\mathcaligr{B}(S^{(j+k)})$, we find that
%
\begin{eqnarray}\label{eqq5}
&& \bigl\vert\bigl[ \Pi_{n}^{(k,m),(j)}(\mu)-\Pi_{n}^{(k,m),(j)}(\nu
) \bigr]
(f) \bigr\vert
\nonumber\\
&&\qquad\leq\int \bigl\vert\bigl[ \bigl( \overline{\Phi}%
{}^{(j+(k-1),j+1)}_{n}\bigl(\mu^{(j)}\bigr) \bigr)\\
&&\qquad\quad\hspace*{14.51pt}{} - \bigl( \overline{\Phi}%
{}^{(j+(k-1),j+1)}_{n}\bigl(\nu^{(j)}\bigr) \bigr) \bigr] (g) \bigr\vert
\Gamma
_{j+k}(f,dg).\nonumber
\end{eqnarray}
We let $\mathcaligr{F}_{n}^{m,j}$ be the sigma field given by
\[
\mathcaligr{F}_{n}^{m,j}=\sigma\bigl( X_{p}^{(l)} \dvtx0\leq p\leq
n, 0\leq l\leq
m, l\not=j \bigr).
\]
Combining the generalized Minkowski integral inequality presented in
Lemma %
\ref{mink} with the inequality (\ref{Laref}), we prove that
\begin{eqnarray*}
&&\mathbb{E} \bigl( \bigl\vert\bigl[ \Pi_{n}^{(k,m),(j)}\bigl(\bigl(\overline
{\eta}%
^{(l)}\bigr)_{l}\bigr)-\Pi_{n}^{(k,m),(j)} \bigl( \bigl( \overline{\Phi}{}^{(l)}\bigl(%
\overline{\eta}^{(l-1)}\bigr) \bigr) _{l} \bigr) \bigr] (f)
\bigr\vert
^{r} \vert\mathcaligr{F}_{n}^{m,j} \bigr) ^{{1/r}}
\\
&&\qquad\leq\int\mathbb{E} \bigl( \bigl\vert\bigl[ \bigl(
\overline{%
\Phi}{}^{(j+(k-1),j+1)}_{n}\bigl(\overline{\eta}^{(j)}\bigr) \bigr)\\
&&\hspace*{63.1pt}{} - \bigl(
\overline{%
\Phi}{}^{(j+(k-1),j+1)}_{n}\bigl(\overline{\Phi}{}^{(j)}\bigl(\overline{\eta}%
{}^{(j-1)}\bigr)\bigr) \bigr) \bigr] (g) \bigr\vert^{r} \vert\mathcaligr
{F}%
_{n}^{m,j} \bigr) ^{{1}/{r}}
\times\Gamma_{j+k}(f,dg)
\\
&&\qquad\leq\frac{e(r)}{\sqrt{n+1}} A B^{k} \Vert f\Vert.%
\end{eqnarray*}
Notice that the decomposition (\ref{Pi}) can be rewritten for any
$f\in%
\mathcaligr{B} ( \prod_{l=k}^{m}S^{(l)} ) $ in the following form:
%
\begin{eqnarray}\label{decomp2}
&& \bigl[ \Pi_{n}^{(k,m)}(\mu)-\Pi_{n}^{(k,m)}(\nu) \bigr] (f)
\nonumber\\[-8pt]\\[-8pt]
&&\qquad= \sum_{j=0}^{m-k} \bigl[ \Pi_{n}^{(k,m),(j)}(\mu)-\Pi
_{n}^{(k,m),(j)}(\nu) \bigr] \bigl( R_{n}^{(k,m),(j)}(\mu,\nu
)(f) \bigr)\nonumber
\end{eqnarray}
with the integral operators $R_{n}^{(k,m),(j)}(\mu,\nu) \dvtx\mathcaligr
{B} (
\prod_{l=k}^{m}S^{(l)} ) \mapsto\mathcaligr{B}(S^{(j+k)})$ given below
\begin{eqnarray*}
&&R_{n}^{(k,m),(j)}(\mu,\nu)(f)(x_{k+j})
\\
&&\qquad= \int f\bigl(x_{k},%
\ldots,x_{k+(j-1)},x_{k+j},x_{k+j+1},\ldots,x_{m}\bigr)
\\
&&\qquad\quad\hspace*{8.43pt}{} \times\Biggl( \prod_{i=0}^{j-1}\Pi
_{n}^{(k,m),(i)}(\nu) \Biggr) (dx_{i+k})\times\Biggl(
\prod_{i=j+1}^{m-k}\Pi_{n}^{(k,m),(i)}(\mu)(dx_{i+k}) \Biggr).%
\end{eqnarray*}
Using the fact that the pair of measures
\[
\bigotimes_{i=0}^{j-1}\Pi_{n}^{(k,m),(i)} \bigl( \bigl( \overline
{\Phi}%
^{(l)}\bigl(\eta^{(l-1)}\bigr) \bigr) _{l} \bigr) \quad\mbox{and}\quad
\bigotimes_{i=j+1}^{m-k}\Pi_{n}^{(k,m),(i)}\bigl(\bigl(\eta^{(l)}\bigr)_{l}\bigr)
\]
only depend on the distribution flow $ ( \overline{\Phi}%
{}^{(i)}(\eta^{(i-1)}) ) _{0\leq i\leq j-1}$ and $(\eta
^{(i)})_{j+1\leq
i\leq m-k}$, we find that the random functions
\[
f_{n}^{(k,m),(j)}:=R_{n}^{(k,m),(j)} \bigl( \bigl(\overline{\eta}%
^{(l)}\bigr)_{l}, \bigl( \overline{\Phi}{}^{(l)}\bigl(\overline{\eta
}^{(l-1)}\bigr) \bigr)
_{l} \bigr) (f)\in\mathcaligr{B}\bigl(S^{(j+k)}\bigr)
\]
do not depend on the distribution flows $\eta^{(j)}$ and $\eta^{(j-1)}$.
This shows that $f_{n}^{(k,m),(j)}$ are measurable with respect to
$\mathcaligr{%
F}_{n}^{m,j}$. From previous calculations (and again using the generalized
Minkowski integral inequality presented in Lemma \ref{mink}) we find that
\begin{eqnarray*}
&&\mathbb{E} \bigl( \bigl\vert\bigl[ \Pi_{n}^{(k,m),(j)}\bigl(\bigl(\overline
{\eta}%
^{(l)}\bigr)_{l}\bigr)-\Pi_{n}^{(k,m),(j)} \bigl( \bigl( \overline{\Phi}{}^{(l)}\bigl(%
\overline{\eta}^{(l-1)}\bigr) \bigr) _{l} \bigr) \bigr]
\bigl(f_{n}^{(k,m),(j)}\bigr)%
\bigr\vert^{r} \vert\mathcaligr{F}_{n}^{m,j} \bigr)
^{{1}/{r}}
\\
&&\qquad\leq\int\Gamma_{j+k}\bigl(f_{n}^{(k,m),(j)},dg\bigr) \\
&&\hspace*{41pt}{}\times \mathbb{E} \bigl( \bigl\vert\bigl[ \bigl( \overline{\Phi
}%
{}^{(j+(k-1),j+1)}_{n}\bigl(\overline{\eta}^{(j)}\bigr) \bigr)\\
&&\qquad\quad\hspace*{37.92pt}{} - \bigl(
\overline{\Phi}%
{}^{(j+(k-1),j+1)}_{n}\bigl(\overline{\Phi}{}^{(j)}\bigl(\overline{\eta
}^{(j-1)}\bigr)\bigr) \bigr) %
\bigr] (g) \bigr\vert^{r} \vert\mathcaligr{F}_{n}^{m,j}\bigr)
^{{1}/{r}} \\
&&\qquad
\leq\frac{e(r)}{\sqrt{n+1}} A B^{k} \Vert f\Vert.%
\end{eqnarray*}
We conclude that for any $f\in\mathcaligr{B}(\prod_{k\leq j\leq m}S^{(j)})$
\begin{eqnarray*}
&&\mathbb{E} \bigl( \bigl\vert\bigl[ \Pi_{n}^{(k,m)}\bigl(\bigl(\overline
{\eta}%
^{(l)}\bigr)_{l}\bigr)-\Pi_{n}^{(k,m)} \bigl( \bigl( \overline{\Phi}{}^{(l)}\bigl(\overline{%
\eta}^{(l-1)}\bigr) \bigr) _{l} \bigr) \bigr] (f) \bigr\vert
^{r} \bigr) ^{{1}/{r}}
\\
&&\qquad\leq(m-k+1) \frac{e(r)}{\sqrt{n+1}} A B^{k} \Vert f\Vert.%
\end{eqnarray*}
Using (\ref{decomp2}), it is now easily checked that for every $f\in
\mathcaligr{B}(E)$
\[
\mathbb{E} \bigl( \vert[ \overline{\omega
}_{n}^{k}(\overline{\eta}%
)-\overline{\omega}_{n}^{k+1}(\overline{\eta}) ] (f)
\vert
^{r} \bigr) ^{{1}/{r}}
\leq(m-k+1) \frac{e(r)}{\sqrt{n+1}} A B^{k} \Vert f\Vert.
\]
Finally, by (\ref{decomp1}) we conclude that
\[
\mathbb{E} \bigl( \bigl\vert\bigl[ \overline{\eta}_{n}-\overline
{\pi}^{[m]}%
\bigr] (f) \bigr\vert^{r} \bigr) ^{{1}/{r}}
\leq\frac{e(r)}{\sqrt{n+1}} A \Vert f\Vert\sum_{k=0}^{m}(m-k+1) B^{k}.
\]
This ends the proof of the theorem.
\end{pf}

\subsection{Concentration analysis}
\label{concentrat}

This section is mainly concerned with exponential bounds
for the deviations of the occupation measures $\overline{\eta}_{n}^{[m]}$
around the limiting tensor product measure $\overline{\pi}{}^{[m]}$. We
restrict our attention to models satisfying the Lipschitz type
condition (%
\ref{lipcn}) for some kernel $\Gamma_{k}$ with uniformly finite support
\[
\sup_{f\in\mathcaligr{B}(S^{(k)})}\operatorname{Card} (\operatorname
{Supp}(\Gamma
_{k}(f,\cdot)) ) <\infty.
\]
To simplify the presentation, we fix a parameter $m\geq1$, and
sometimes we
write $\overline{\eta}_{n}$ instead of $\overline{\eta}_{n}^{[m]}$.
We shall
also use the letters $c_{i}$, $i\geq1$ to denote some finite constants whose
values may vary from line to line but do not depend on the time
parameter $n$.

The main result of this section is the following concentration theorem:
\begin{theorem}
\label{conctheo} There exists a finite constant $\overline{\sigma}%
_{m}<\infty$ such that for any $f\in\mathcaligr{B}_{1}(E_{m})$ and $t>0$
\[
\limsup_{n\rightarrow\infty}\frac{1}{n}\log\mathbb{P} \bigl(
\bigl\vert%
\bigl[ \overline{\eta}_{n}^{[m]}-\overline{\pi}{}^{[m]} \bigr]
(f) \bigr\vert>t \bigr) <-\frac{t^{2}}{2\overline{\sigma}_{m}^{2}}.
\]
\end{theorem}

The proof of this theorem is based on two technical lemmas.
\begin{lem}
\label{leminter} We let $M=(M_{n})_{n\geq1}$ be a random process such that
the following exponential inequality is satisfied for some positive
constants $a,b>0$ and for any $t\geq0$ and $n\geq1$
\[
\mathbb{P}\bigl( \vert M_{n} \vert\geq t \sqrt{n}\bigr)\leq a e^{-{bt^{2}}}.
\]
We consider the collection of random processes $\overline{M}{}^{(k)}=(%
\overline{M}{}^{(k)}_{n})_{n\geq1}$ defined for any $n\geq0$ and
$k\geq0$
by the following formula:
\[
\overline{M}{}^{(k)}_{n+1}:=(n+1) \int\Sigma^{k}(n,dp) \frac{1}{p+1}
M_{p+1},
\]
where $\Sigma^{k}$ is the semigroup associated to the operator $\Sigma$
defined in (\ref{defsigma}). For every $k\geq0$, $n\geq1$, and
$t\geq0$
we have the exponential inequalities:
\[
\mathbb{P} \bigl( \bigl\vert\overline{M}{}^{(k)}_{n} \bigr\vert\geq
t \sqrt{n}%
\bigr) \leq a n^{k} e^{-{bt^{2}/2^{2k}}}.
\]
\end{lem}
\begin{pf}
We prove the lemma by induction on the parameter $k$. For $k=0$, we
have $%
\overline{M}{}^{(0)}_{n+1}:=M_{n+1}$ so that the exponential estimate holds
true with $a(0)=a$ and $b(0)=b$. Suppose we have proved the result at
rank $%
k $. Using the fact that
\begin{eqnarray*}
\overline{M}_{n+1}^{(k+1)} &=& (n+1) \int\Sigma^{k+1}(n,dp) \frac
{1}{p+1}%
M_{p+1} \\
&=& (n+1) \int\Sigma(n,dp) \frac{1}{p+1} \biggl( (p+1)\int\Sigma
^{k}(p,dq) %
\frac{1}{q+1} M_{q+1} \biggr)
\end{eqnarray*}
we prove the recursion formula
\[
\overline{M}{}^{(k+1)}_{n+1}=(n+1) \int\Sigma(n,dp) \frac
{1}{p+1} \overline{M%
}{}^{(k)}_{p+1}.
\]
On the other hand, we have
\[
\frac{1}{2} \frac{\overline{M}{}^{(k+1)}_{n+1}}{\sqrt{n+1}}=\frac
{1}{2} \sqrt{%
n+1} \int\Sigma(n,dp) \frac{1}{\sqrt{p+1}} \frac{\overline
{M}{}^{(k)}_{p+1}}{%
\sqrt{p+1}}
\]
and
\begin{eqnarray*}
\frac{1}{2}\sqrt{n+1} \int\Sigma(n,dp) \frac{1}{\sqrt{p+1}}
&=& \frac{1}{2%
\sqrt{n+1}}\sum_{p=0}^{n} \frac{1}{\sqrt{p+1}} \\
&\leq& \frac{1}{2\sqrt{n+1}}\sum_{p=0}^{n} \int_{p}^{p+1}\frac
{1}{\sqrt{t}}\,
dt=1.
\end{eqnarray*}
Under the induction hypothesis, we have for any $0\leq p\leq n$
\[
\mathbb{P} \bigl( \bigl\vert\overline{M}{}^{(k)}_{p+1} \bigr\vert
\geq t \sqrt{%
p+1} \bigr) \leq a (n+1)^{k} e^{-{bt^{2}/2^{2k}}}.
\]
This implies that
\begin{eqnarray*}
\mathbb{P} \biggl( \frac{1}{2} \frac{\overline
{M}{}^{(k+1)}_{n+1}}{\sqrt{n+1}}%
>t \biggr)
&\leq&
\mathbb{P} \bigl( \exists0\leq p\leq n \dvtx\overline{M}%
{}^{(k)}_{p+1}>t\sqrt{p+1} \bigr) \\
&\leq&
a (n+1) (n+1)^{k} e^{-{bt^{2}/2^{2k}}}
\end{eqnarray*}
from which we conclude that
\[
\mathbb{P} \bigl( \overline{M}{}^{(k+1)}_{n+1}>t \sqrt{n+1} \bigr)
\leq
a (n+1)^{k+1} e^{-{bt^{2}/2^{2(k+1)}}}.
\]
This ends the proof of the lemma.
\end{pf}
\begin{lem}
\label{lemF} For every $l_{1}<l_{2}$, there exists some nonincreasing
function
\[
N \dvtx t\in[0,\infty)\quad\mapsto\quad N(t)\in[0,\infty)
\]
such that for every $n\geq N(t)$ and any function $f\in\mathcaligr{B}%
_{1}(S^{(l_{2})})$ we have
\begin{eqnarray*}
&&\mathbb{P} \bigl( \sqrt{n+1} \bigl\vert\bigl[ \overline{\Phi}%
{}^{(l_{2},l_{1}+1)}_{n}\bigl(\overline{\eta}^{(l_{1})}\bigr)-\overline{\Phi}%
{}^{(l_{2},l_{1})}_{n}\bigl(\overline{\Phi}{}^{(l_{1})}\bigl(\overline{\eta
}^{(l_{1}-1)}\bigr)\bigr)%
\bigr] (f) \bigr\vert>t \bigr) \\
&&\qquad
\leq\bigl(c_{1}(n+1)\bigr)^{(l_{2}-l_{1})} \exp{ (
-{c_{2}t^{2}/c_{3}^{l_{2}-l_{1}}%
} ) }.
\end{eqnarray*}
\end{lem}

Before getting into the details of the proof of this lemma, it is
interesting to mention a direct consequence of the above exponential
estimates. First, we observe that $N(t\sqrt{n+1})\leq N(t)$ so that
for any $%
t>0$ and $n\geq N(t)$ we have
\begin{eqnarray*}
&&\mathbb{P} \bigl( \bigl\vert\bigl[ \overline{\Phi
}{}^{(l_{2},l_{1}+1)}_{n}\bigl(%
\overline{\eta}{}^{(l_{1})}\bigr)-\overline{\Phi
}{}^{(l_{2},l_{1})}_{n}\bigl(\overline{%
\Phi}{}^{(l_{1})}\bigl(\overline{\eta}^{(l_{1}-1)}\bigr)\bigr) \bigr] (f)
\bigr\vert
>t \bigr) \\
&&\qquad
\leq\bigl(c_{1}(n+1)\bigr)^{(l_{2}-l_{1})} \exp \bigl( -{%
c_{2}(n+1)t^{2}/c_{3}^{l_{2}-l_{1}}} \bigr) .%
\end{eqnarray*}
Using the decomposition
\[
\eta_{n}^{(k)}-\pi^{(k)}=\sum_{l=0}^{k} \bigl[ \overline{\Phi}%
{}^{(k,l+1)}_{n}\bigl(\eta^{(l)}\bigr)-\overline{\Phi}{}^{(k,l+1)}_{n}\bigl(\overline
{\Phi}%
{}^{(l)}\bigl(\eta^{(l-1)}\bigr)\bigr) \bigr]
\]
we prove the following inclusion of events:
\begin{eqnarray*}
&&\hspace*{-3pt} \bigl\{ \bigl\vert\bigl[\overline{\eta}_{n}^{(k)}-\pi
^{(k)}\bigr](f) \bigr\vert>t \bigr\}
\\
&&\hspace*{-3pt}\qquad{}\subset\bigl\{ \exists0\leq l\leq k \dvtx \bigl\vert\bigl[ \overline
{\Phi}%
{}^{(k,l+1)}_{n}\bigl(\eta^{(l)}\bigr)-\overline{\Phi}_{n}{}^{(k,l+1)}\bigl(\overline
{\Phi}%
{}^{(l)}\bigl(\eta^{(l-1)}\bigr)\bigr) \bigr] (f) \bigr\vert>t/(k+1)
\bigr\}.
\end{eqnarray*}
By Lemma \ref{lemF} we can find a sufficiently large integer $N(t)$
that may
depend on the parameter $k$ and such that for every $n\geq N(t)$
\begin{eqnarray*}
&&\mathbb{P} \bigl( \bigl\vert\bigl[\overline{\eta}_{n}^{(k)}-\pi
^{(k)}\bigr](f) \bigr\vert>t \bigr) \\
&&\qquad
\leq\sum_{0\leq l\leq k} \mathbb{P} \biggl( \bigl\vert\bigl[
\overline{\Phi}%
{}^{(k,l+1)}_{n}\bigl(\overline{\eta}^{(l)}\bigr)-\overline{\Phi}{}^{(k,l)}_{n}\bigl(%
\overline{\Phi}{}^{(l)}\bigl(\overline{\eta}^{(l-1)}\bigr)\bigr) \bigr] (f)
\bigr\vert>%
\frac{t}{k+1} \biggr) \\
&&\qquad
\leq(k+1) \bigl(c_{1}(n+1)\bigr)^{k} e^{-{(n+1)t^{2}c_{2}/((k+1)^{2}c_{3}^{k}})}.
\end{eqnarray*}
This clearly implies the existence of some finite constant $\sigma
_{k}<\infty$ such that
\[
\limsup_{n\rightarrow\infty}\frac{1}{n}\log{\mathbb{P} \bigl(
\bigl\vert\bigl[%
\overline{\eta}_{n}^{(k)}-\pi^{(k)}\bigr](f) \bigr\vert>t \bigr)
}<-\frac{t^{2}%
}{2\sigma_{k}^{2}}.
\]
\begin{pf*}{Proof of Lemma \protect\ref{lemF}}
Using Lemma \ref{phibar}, we find that
\begin{eqnarray*}
&& \bigl\vert\bigl[ \overline{\Phi}{}^{(l_{2},l_{1}+1)}_{n}\bigl(\overline
{\eta}%
^{(l_{1})}\bigr)-\overline{\Phi}{}^{(l_{2},l_{1})}_{n}\bigl(\overline{\Phi}{}^{(l_{1})}\bigl(%
\overline{\eta}^{(l_{1}-1)}\bigr)\bigr) \bigr] (f) \bigr\vert\\
&&\qquad
\leq\int \bigl\vert\bigl[ \overline{\eta
}_{p}^{(l_{1})}-%
\overline{\Phi}{}^{(l_{1})}_{p}\bigl(\overline{\eta}^{(l_{1}-1)}\bigr) \bigr]
(g) \bigr\vert\overline{\Gamma}{}^{(l_{2},l_{1}+1)}((n,f),d(p,g)).%
\end{eqnarray*}
Arguing as in (\ref{chernov}), we find that for any $g\in\mathcaligr{B}%
(S^{(l_{1})})$, we have
%
\begin{equation}\label{eqq1b}
\bigl\vert\bigl[ \overline{\eta}_{p}^{(l_{1})}-\overline{\Phi}%
{}^{(l_{1})}_{p}\bigl(\overline{\eta}^{(l_{1}-1)}\bigr) \bigr] (g) \bigr\vert
\leq\frac{%
\vert M_{p+1}^{(l_{1})}(g) \vert}{p+1}+c_{1} \frac{\log
{(p+2)}}{%
p+2} \Vert g\Vert
\end{equation}
with a sub-Gaussian process $M_{n}^{(l_{1})}(g)$ satisfying the following
exponential inequality for any $t>0$ and any time parameter $n\geq1$:
\[
\mathbb{P}\bigl( \bigl\vert M_{n}^{(l_{1})}(g) \bigr\vert\geq t\sqrt
{n}\bigr)\leq
2 \exp{ ( -c_{2}t^{2}/\Vert g\Vert^{2} ) }.
\]
We notice that
\begin{eqnarray*}
\frac{1}{n+2}\sum_{p=0}^{n}\frac{(\log{(p+2)})^{k}}{p+2}
&\leq&
\frac{(\log{%
(n+2)})^{k}}{n+2}\sum_{p=0}^{n}\frac{1}{p+2} \\
&\leq&
\frac{(\log{(n+2)})^{k}}{n+2}\sum_{p=0}^{n}\int
_{p+1}^{p+2}\frac{1}{t}%
\,dt\\
&=&\frac{(\log{(n+2)})^{k+1}}{n+2}.
\end{eqnarray*}
This implies that
\[
\int\Sigma(n,dp) \frac{\log{(p+2)}}{p+2}\leq2 \frac{(\log
{(n+2)})^{2}}{n+2}.
\]
More generally for any $k\geq0$, we have that
\[
\int\Sigma^{k}(n,dp) \frac{\log{(p+2)}}{p+2}\leq2^{k} \frac
{(\log{(n+2)}%
)^{k+1}}{n+2}
\]
from which we prove that
%
\begin{eqnarray}\label{eqq2b}
&& \int\frac{\log{(p+2)}}{p+2} \Vert g\Vert\overline
{\Gamma}%
{}^{(l_{2},l_{1}+1)}((n,f),d(p,g)) \nonumber\\
&&\qquad
\leq2^{(l_{2}-l_{1})} \frac{(\log
{(n+2)})^{(l_{2}-l_{1})+1}}{n+2%
} \int\Vert g\Vert\Gamma_{l_{2},l_{1}+1}(f,dg) \nonumber\\[-8pt]\\[-8pt]
&&\qquad
\leq2^{(l_{2}-l_{1})} \frac{(\log
{(n+2)})^{(l_{2}-l_{1})+1}}{n+2%
} \biggl( \prod_{l_{1}<i\leq l_{2}}\Lambda_{i} \biggr)\nonumber\\
&&\qquad \leq
c_{3}^{(l_{2}-l_{1})} \frac{(\log
{(n+2)})^{(l_{2}-l_{1})+1}}{n+2%
}.\nonumber
\end{eqnarray}
For any $g\in\mathcaligr{B}(S^{(l_{1})})$ we set
\[
\overline{\mathcaligr{M}}{}^{(l_{1},l_{2})}_{n+1}(g):=\int%
\Sigma^{(l_{2}-l_{1})}(n,dp) \frac{ \vert
M_{p+1}^{(l_{1})}(g) \vert}{p+1}.
\]
Using Lemma \ref{leminter}, we prove that
\[
\mathbb{P} \bigl( \overline{\mathcaligr
{M}}{}^{(l_{1},l_{2})}_{n+1}(g)>t \bigr)
\leq2 (n+1)^{(l_{2}-l_{1})} \exp \bigl( -%
c_{2}(n+1)t^{2}/\bigl[2^{2(l_{2}-l_{1})}\Vert g\Vert^{2}\bigr] \bigr).
\]
We observe that
\[
\int\frac{1}{p+1} \bigl\vert
M_{p+1}^{(l_{1})}(g) \bigr\vert%
\overline{\Gamma}{}^{(l_{2},l_{1}+1)}((n,f),d(p,g))= \int
\overline{\mathcaligr{M}}{}^{(l_{1},l_{2})}_{n+1}(g) \Gamma_{l_{2},l_{1}+1}(f,dg).
\]
In addition, using (\ref{eqq1b}) and (\ref{eqq2b}) we find that
%
\begin{eqnarray}\label{eqq3}
&&
\bigl\vert\bigl[ \overline{\Phi}{}^{(l_{2},l_{1}+1)}_{n}\bigl(\overline
{\eta}%
{}^{(l_{1})}\bigr)-\overline{\Phi}_{n}{}^{(l_{2},l_{1})}\bigl(\overline{\Phi}{}^{(l_{1})}\bigl(%
\overline{\eta}^{(l_{1}-1)}\bigr)\bigr) \bigr] (f) \bigr\vert\nonumber\\[-8pt]\\[-8pt]
&&\qquad
\leq\int\overline{\mathcaligr
{M}}{}^{(l_{1},l_{2})}_{n+1}(g) %
\Gamma_{l_{2},l_{1}+1}(f,dg)+\varepsilon_{l_{1},l_{2}}(n)\nonumber
\end{eqnarray}
with
\[
\varepsilon_{l_{1},l_{2}}(n):=c_{1}c_{3}^{(l_{2}-l_{1})}
\frac{%
(\log{(n+2)})^{(l_{2}-l_{1})+1}}{n+2}.
\]
Using the inclusion of events
\begin{eqnarray*}
&& \biggl\{ \int\overline{\mathcaligr{M}}{}^{(l_{1},l_{2})}_{n+1}(g) \Gamma
_{l_{2},l_{1}+1}(f,dg)>t \biggr\} \\
&&\qquad{}
\subset\bigl\{ \exists g\in\operatorname{Supp}(\Gamma
_{l_{2},l_{1}+1}(f,%
\cdot)) \mbox{ such that } \overline{\mathcaligr{M}}%
{}^{(l_{1},l_{2})}_{n+1}(g)>t\Vert g\Vert/ (
\Lambda_{l_{2},l_{1}+1} ) \bigr\}%
\end{eqnarray*}
we find that
\begin{eqnarray*}
&&\mathbb{P} \biggl( \int\overline{\mathcaligr
{M}}{}^{(l_{1},l_{2})}_{n+1}(g) %
\Gamma_{l_{2},l_{1}+1}(f,dg)>t \biggr) \\
&&\qquad
\leq S_{l_{2},l_{1}+1}(f) \mathbb{P} \bigl( \overline{\mathcaligr{M}}%
{}^{(l_{1},l_{2})}_{n+1}(g)>t\Vert g\Vert/ (
\Lambda_{l_{2},l_{1}+1} ) \bigr).%
\end{eqnarray*}
Finally, under our assumptions we have
\begin{eqnarray*}
S_{l_{2},l_{1}+1}(f) &=& \operatorname{Card} ( \operatorname{Supp}(\Gamma
_{l_{2},l_{1}+1}(f,\cdot)) ) \\
&\leq& \prod_{l_{1}+1\leq k\leq l_{2}}\sup_{f\in\mathcaligr
{B}(S^{(k)})}%
\operatorname{Card} ( \operatorname{Supp}(\Gamma_{k}(f,\cdot)) )
\leq
c_{4}^{(l_{2}-l_{1})}
\end{eqnarray*}
from which we check that
\begin{eqnarray*}
&&\mathbb{P} \biggl( \int\overline{\mathcaligr
{M}}{}^{(l_{1},l_{2})}_{n+1}(g) %
\Gamma_{l_{2},l_{1}+1}(f,dg)>t \biggr) \\
&&\qquad
\leq\bigl(c_{5}(n+1)\bigr)^{(l_{2}-l_{1})} \exp \bigl( -{%
c_{6}(n+1)t^{2}/c_{7}^{(l_{2}-l_{1})}} \bigr) .%
\end{eqnarray*}
Using (\ref{eqq3}), we conclude that
\begin{eqnarray*}
&&\mathbb{P} \bigl( \bigl\vert\bigl[ \overline{\Phi
}{}^{(l_{2},l_{1}+1)}_{n}\bigl(%
\overline{\eta}^{(l_{1})}\bigr)-\overline{\Phi
}{}^{(l_{2},l_{1})}_{n}\bigl(\overline{\Phi%
}{}^{(l_{1})}\bigl(\overline{\eta}^{(l_{1}-1)}\bigr)\bigr) \bigr] (f) \bigr\vert
>t+\varepsilon_{l_{1},l_{2}}(n) \bigr) \\
&&\qquad
\leq\bigl(c_{5}(n+1)\bigr){}^{(l_{2}-l_{1})} \exp \bigl( -{%
c_{6}(n+1)t^{2}/c_{7}^{(l_{2}-l_{1})}} \bigr) .%
\end{eqnarray*}
To take the final step, we observe that
\begin{eqnarray*}
&&\mathbb{P} \bigl( \sqrt{n+1} \bigl\vert\bigl[ \overline{\Phi}%
{}^{(l_{2},l_{1}+1)}_{n}\bigl(\overline{\eta}^{(l_{1})}\bigr)-\overline{\Phi}%
{}^{(l_{2},l_{1})}_{n}\bigl(\overline{\Phi}{}^{(l_{1})}\bigl(\overline{\eta
}^{(l_{1}-1)}\bigr)\bigr)%
\bigr] (f) \bigr\vert
>t+\sqrt{n+1} \varepsilon_{l_{1},l_{2}}(n) \bigr) \\
&&\qquad
\leq\mathbb{P} \biggl( \bigl\vert\bigl[ \overline{\Phi
}{}^{(l_{2},l_{1}+1)}_{n}\bigl(%
\overline{\eta}^{(l_{1})}\bigr)-\overline{\Phi
}{}^{(l_{2},l_{1})}_{n}\bigl(\overline{\Phi%
}{}^{(l_{1})}\bigl(\overline{\eta}^{(l_{1}-1)}\bigr)\bigr) \bigr] (f) \bigr\vert
\\[-4pt]
&&\qquad\quad\hspace*{128pt}
>\frac{t}{%
\sqrt{n+1}}+\varepsilon_{l_{1},l_{2}}(n) \biggr).%
\end{eqnarray*}
We also notice that for any $t>0$ we can find some nonincreasing
function $%
N(t)$ such that
\[
\forall n\geq N(t)\qquad \sqrt{n+1} \varepsilon_{l_{1},l_{2}}(n)<t.
\]
This implies that for any $n\geq N(t)$ we have
\begin{eqnarray*}
&&\mathbb{P} \bigl( \sqrt{n+1} \bigl\vert\bigl[ \overline{\Phi}%
{}^{(l_{2},l_{1}+1)}_{n}\bigl(\overline{\eta}^{(l_{1})}\bigr)-\overline{\Phi}%
{}^{(l_{2},l_{1})}_{n}\bigl(\overline{\Phi}{}^{(l_{1})}\bigl(\overline{\eta
}^{(l_{1}-1)}\bigr)\bigr)%
\bigr] (f) \bigr\vert>2t \bigr) \\
&&\qquad\leq\bigl(c_{5}(n+1)\bigr)^{(l_{2}-l_{1})} \exp \bigl( -{%
c_{6}t^{2}/c_{7}^{(l_{2}-l_{1})}} \bigr).
\end{eqnarray*}
The end of the proof is now straightforward.
\end{pf*}

We are now in position to prove Theorem \ref{conctheo}.
\begin{pf*}{Proof of Theorem \protect\ref{conctheo}}
We use the same notation as we used in the proof of Theorem \ref{theo411}.
Using (\ref{eqq5}) we find that
\begin{eqnarray*}
&& \bigl\vert\bigl[ \Pi_{n}^{(k,m),(j)}(\mu)-\Pi_{n}^{(k,m),(j)}(\nu
) \bigr]
(f) \bigr\vert>t \\
&&
\quad\Longrightarrow\quad\exists g\in\operatorname{Supp}(\Gamma_{j+k}(f,\cdot
)) \dvtx
\bigl\vert\bigl[ \bigl( \overline{\Phi
}{}^{(j+(k-1),j+1)}_{n}\bigl(%
\mu^{(j)}\bigr) \bigr)\\
&&\qquad\quad\hspace*{121pt}{} - \bigl( \overline{\Phi}{}^{(j+(k-1),j+1)}_{n}\bigl(\nu
^{(j)}\bigr) \bigr) \bigr] (g) \bigr\vert>t\Vert g\Vert/\Lambda_{j+k}.
\end{eqnarray*}
Therefore, using Lemma \ref{lemF} we can find a nonincreasing function
$%
N(t) $ (that may depend on the parameter $k$), such that for every
$n\geq
N(t)$ and any $f\in\mathcaligr{B}_{1}(S^{(j+k)})$ we have
\begin{eqnarray*}
&&\mathbb{P} \bigl( \sqrt{n+1} \bigl\vert\bigl[ \Pi
{}^{(k,m),(j)}_{n}(\mu
)-\Pi_{n}^{(k,m),(j)}(\nu) \bigr] (f) \bigr\vert>t \bigr) \\
&&\qquad
\leq\bigl(c_{1}(n+1)\bigr)^{(k-1)} \exp \bigl(
-{c_{2}t^{2}/c_{3}^{(k-1)}} \bigr) .%
\end{eqnarray*}
In much the same way, by the decomposition (\ref{decomp2}) we find the
following assertion:
\begin{eqnarray*}
&&
\bigl\vert\bigl[ \Pi_{n}^{(k,m)}(\mu)-\Pi_{n}^{(k,m)}(\nu) \bigr]
(f) \bigr\vert>t \\
&&\quad
\Longrightarrow\quad\exists0\leq j\leq(m-k) \dvtx
\bigl\vert\bigl[ \Pi_{n}^{(k,m),(j)}(\mu) -\Pi
_{n}^{(k,m),(j)}(\nu)%
\bigr]\\
&&\qquad\quad\hspace*{126pt}{}\times \bigl( R_{n}^{(k,m),(j)}(\mu,\nu)(f) \bigr) \bigr\vert
>t/(m-k+1).%
\end{eqnarray*}
Since $R_{n}^{(k,m),(j)}(\mu,\nu)$ maps $\mathcaligr{B}_{1}(%
\prod_{l=k}^{m}S^{(l)})$ into $\mathcaligr{B}_{1}(S^{(j+k)})$ we have
for every
parameter $n\geq N(t)$
\begin{eqnarray*}
&&
\mathbb{P} \bigl( \sqrt{n+1} \bigl\vert\bigl[ \Pi
_{n}^{(k,m)}\bigl(\bigl(\overline{%
\eta}^{(l)}\bigr)_{l}\bigr)-\Pi_{n}^{(k,m)}\bigl(\bigl(\overline{\Phi}{}^{(l)}\bigl(\overline
{\eta}%
^{(l-1)}\bigr)\bigr)_{l}\bigr) \bigr] (f) \bigr\vert>t \bigr) \\
&&\qquad
\leq(m-k+1) \bigl(c_{1}(n+1)\bigr)^{k-1} \exp \bigl( -
c_{2}t^{2}/\bigl((m-k+1)^{2}c_{3}^{k-1}\bigr) \bigr) .%
\end{eqnarray*}
In summary, we have proved that there exists some nonincreasing
function $%
N(t)$ that may depend on the parameter $m$ such that for any $0\leq
k\leq m$%
, any $f\in\mathcaligr{B}_{1}(E)$, and any $n\geq N(t)$ we have
\begin{eqnarray*}
&&
\mathbb{P} \bigl( \sqrt{n+1} \bigl\vert\bigl[ \overline{\pi}%
^{[k-1]}\otimes\bigl\{ \Pi_{n}^{(k,m)}\bigl(\bigl(\overline{\eta
}^{(l)}\bigr)_{l}\bigr)-%
\Pi_{n}^{(k,m)}\bigl(\bigl(\overline{\Phi}{}^{(l)}\bigl(\overline{\eta
}^{(l-1)}\bigr)\bigr)_{l}\bigr) %
\bigr\} \bigr] (f) \bigr\vert>t \bigr) \\
&&\qquad
\leq\bigl(c_{4}(n+1)\bigr)^{m} \exp ( -c_{5}t^{2}/c_{6}^{m}) .
\end{eqnarray*}
Let $(U_{n})_{n\geq1}$ be a collection of $[0,1]$-valued random variables
such that for any $t$ there exists some nonincreasing function $N(t)$, so
that for $n\geq N(t)$
\[
\mathbb{P}\bigl(\sqrt{n} U_{n}\geq t\bigr)\leq a n^{\alpha} e^{-t^{2}b}
\]
for some integer $\alpha\geq1$ and some pair of positive constants $(a,b)$.
In this situation, we can find a nonincreasing function $N^{\prime
}(t)$ and
a pair of positive constants $(a^{\prime},b^{\prime})$ such that
\[
\forall n\geq N^{\prime}(t) \mathbb{P}\qquad \Biggl( \sum
_{p=1}^{n}U_{p}>\sqrt{n%
} t \Biggr) \leq a^{\prime} n^{\alpha+1} e^{-t^{2}b^{\prime}}.
\]
To prove this claim, we simply use the fact that for any $n\geq N(t)$ we
have
\[
\frac{1}{\sqrt{n}}\sum_{p=1}^{n}U_{p}\leq\frac{N(t)}{\sqrt
{n}}+\frac{1}{%
\sqrt{n}}\sum_{p=N(t)}^{n}\frac{1}{\sqrt{p}} \bigl(\sqrt{p}U_{p}\bigr)
\quad\mbox{and}\quad \frac{1}{2\sqrt{n}}\sum_{p=1}^{n}\frac{1}{\sqrt{p}}\leq1.
\]
This yields that for any $n\geq N(t)$
\[
\mathbb{P} \Biggl( \frac{1}{\sqrt{n}}\sum_{p=1}^{n}U_{p}>t+\frac
{N(t)}{\sqrt{n%
}} \Biggr) \leq\sum_{p=N(t)}^{n}\mathbb{P} \bigl( \sqrt
{p}U_{p}>t/2 \bigr).
\]
We let $N^{\prime}(t)$ be the smallest integer $n$ such that
$N(t)/\sqrt{n}%
\leq t$. Recalling that $N(t)$ is a nondecreasing function, we find that
for any $s\geq t$
\[
N(t)/\sqrt{n}\leq t\quad\Longrightarrow\quad N(s)/\sqrt{n}\leq N(t)/\sqrt
{n}\leq t\leq
s\quad\Longrightarrow\quad N(s)/\sqrt{n}\leq s.
\]
This implies that $N^{\prime}(s)\leq N^{\prime}(t)$. Thus, we have
constructed a nonincreasing function $N^{\prime}(t)$ such that for any
$%
n\geq N^{\prime}(t)$
\[
\mathbb{P} \Biggl( \frac{1}{\sqrt{n}}\sum_{p=1}^{n}U_{p}>2t \Biggr)
\leq
a n^{\alpha+1} e^{-t^{2}b/4}.
\]
This ends the proof of the assertion with $(a^{\prime},b^{%
\prime})=(a,b/2^{4})$. Applying this property to the decomposition
(\ref%
{yenamarre}), we can find a nonincreasing function $N(t)$ such that for any
$n\geq N(t)$ and any $0\leq k\leq m$
\[
\mathbb{P} \bigl( \sqrt{n+1} \vert[ \overline{\omega}%
{}^{k}_{n}(\eta)-\overline{\omega}_{n}^{k+1}(\eta) ] (f)
\vert
>t \bigr) \leq\bigl(c_{7}(n+1)\bigr)^{m+1} \exp (
-{c_{8}t^{2}/c_{9}^{m}} ).
\]
The end of the proof of the theorem is now a direct consequence of the
decomposition (\ref{decomp1}).
\end{pf*}

\section{Feynman--Kac semigroups}
\label{FKS}

In Section \ref{unifth}, we established a uniform convergence theorem under
the assumption that the time averaged semigroup $\overline{\Phi}{}^{(k,l)}$
introduced in Section \ref{timeav} is exponentially stable; that is, it
satisfies (\ref{exponentialstability}). In this section, we study the
mappings $\overline{\Phi}{}^{(k,l)}$ associated with the Feynman--Kac
transformations discussed in (\ref{fksgg}). We provide necessary conditions
ensuring that (\ref{exponentialstability}) is satisfied in this case.

\subsection{Description of the models}

To precisely describe these mappings we need a few definitions.
\begin{defi}
We denote by $\Psi_{l}^{G}$ the Boltzman--Gibbs transformation associated
with a positive potential function $G$ on $S^{(l)}$, and defined for
any $%
f\in\mathcaligr{B}(S^{(l)})$ by the following formula:
\[
\Psi_{l}^{G}(\eta_{p})(f)=\eta_{p}(Gf)/\eta_{p}(G).
\]
We let $Q_{l}$ be the integral operator from $\mathcaligr{B }(S^{(l)})$
into $%
\mathcaligr{B }(S^{(l-1)})$ given by
%
\begin{equation} \label{dfQ}
\forall f\in\mathcaligr{B}\bigl(S^{(l)}\bigr)\qquad Q_{l}(f):=G_{l-1}\times
L_{l}(f)\in%
\mathcaligr{B}\bigl(S^{(l-1)}\bigr).
\end{equation}
\end{defi}

By definition of the mappings $\Phi_{l}$ given in (\ref{fksg}), it is easy
to check that
%
\begin{eqnarray}\label{fksgg}
\overline{\Phi}{}^{(l)}(\eta)=\overline{\Psi}{}^{(l), Q_{l}(1)}(\eta
)L_{l} \hspace*{50pt}\nonumber\\[-8pt]\\[-8pt]
\eqntext{\mbox{with } \forall n\geq0 \overline{\Psi
}{}^{(l),Q_{l}(1)}_{n}(%
\eta)=\dfrac{1}{n+1}\displaystyle\sum_{p=0}^{n}\Psi^{Q_{l}(1)}_{l}(\eta_{p}).}
\end{eqnarray}
\begin{defi}
We let $\overline{\Phi}{}^{(k,l)}$ be the semigroup associated with the
Feyn\-man--Kac transformations $\Phi_{l}$ discussed in (\ref{fksgg}),
and we
denote by
\[
Q_{l,k}=Q_{l}Q_{l+1}\cdots Q_{k}
\]
the semigroup associated with the integral operator $Q_{l}$ introduced
in (%
\ref{dfQ}).
\end{defi}
\begin{prop}
\label{laprop} For any $l\leq k$ we have that
%
\begin{equation}\label{propphi}
\overline{\Phi}{}^{(k,l)}(\eta)=\overline{\Psi}{}^{(k,l)}(\eta
)P_{l,k} \qquad
\mbox{with } P_{l,k}(f)=\frac{Q_{l,k}(f)}{Q_{l,k}(1)},
\end{equation}
and the mapping $\overline{\Psi}{}^{(k,l)}$ from $\mathcaligr
{P}(S^{(l-1)})^{%
\mathbb{N}}$ into itself given below:
\begin{eqnarray*}
\overline{\Psi}{}^{(k,l)} &=& \overline{\Psi}{}^{(l),H_{l,k}}\circ
\overline{\Psi}%
{}^{(k-1,l)} \\
&=& \overline{\Psi}{}^{(l),H_{l,k}}\circ\overline{\Psi
}{}^{(l),H_{l,k-1}}\circ%
\cdots\circ\overline{\Psi}{}^{(l),H_{l,l}} \qquad\mbox{with }
H_{l,k}:=%
\frac{Q_{l,k}(1)}{Q_{l,k-1}(1)}.
\end{eqnarray*}
For $l=k$, we use the conventions $\overline{\Psi
}{}^{(k-1,l)}=\overline{\Psi}%
{}^{(l-1,l)}=\mathrm{Id}$ and $Q_{l,k-1}(1)=Q_{l,l-1}(1)=1$, so that $%
H_{l,l}=Q_{l,l}(1)=Q_{l}(1)$ and $\overline{\Psi}{}^{(l,l)}=\overline
{\Psi}%
{}^{(l),Q_{l}(1)}$.
\end{prop}
\begin{pf}
We prove the proposition by induction on the parameter $m=(k-l)$. For $k=l$,
we clearly have
\[
P_{l,l}(f)=\frac{Q_{l}(f)}{Q_{l}(1)}=L_{l}(f)
\]
and
\[
\overline{\Psi}{}^{(l,l)}=\overline{\Psi
}{}^{(l),Q_{l}(1)}\quad\Longrightarrow\quad
\overline{\Phi}{}^{(l)}(\eta)=\overline{\Psi}{}^{(l,l)}(\eta)P_{l,l}.
\]
Suppose we have proved formula (\ref{propphi}) for some $m=(k-l)\geq
0$. To
check the result at level $m+1=(k-l)+1=((k+1)-l)$, we first observe that
\[
\overline{\Phi}{}^{(k+1)} \bigl( \overline{\Phi}{}^{(k,l)}(\eta)
\bigr) =%
\overline{\Psi}{}^{(k+1),Q_{k+1}(1)}\bigl(\overline{\Phi}{}^{(k,l)}(\eta
)\bigr)P_{k+1,k+1}.
\]
For any $\mu\in\mathcaligr{P}(S^{(k)})$, we also have that
\[
\overline{\Psi}{}^{(k+1),Q_{k+1}(1)}_{n}(\mu)(P_{k+1}(f))=\frac
{1}{n+1}%
\sum_{p=0}^{n}\frac{\mu_{p}(Q_{k+1}(f))}{\mu_{p}(Q_{k+1}(1))}
\]
so that
\[
\overline{\Psi}{}^{(k+1),Q_{k+1}(1)}_{n}\bigl(\overline{\Phi}{}^{(k,l)}(\eta
)\bigr)P_{k+1,k+1}=\frac{1}{n+1}\sum_{p=0}^{n}\frac{\overline{\Phi}%
{}^{(k,l)}_{p}(\eta)(Q_{k+1}(f))}{\overline{\Phi}{}^{(k,l)}_{p}(\eta
)(Q_{k+1}(1))}.
\]
Using the induction hypothesis, we find that
\[
\overline{\Phi}{}^{(k,l)}_{p}(\eta)(Q_{k+1}(f))=\overline{\Psi}%
{}^{(k,l)}(\eta)[P_{l,k}(Q_{k+1}(f))].
\]
We also have
\[
P_{l,k}(Q_{k+1}(f))=\frac{Q_{l,k+1}(1)}{Q_{l,k}(1)}%
P_{l,k+1}(f)=H_{l,k+1} P_{l,k+1}(f)
\]
from which we prove that
\[
\overline{\Psi}{}^{(k,l)}(\eta)[P_{l,k}(Q_{k+1}(f))]=\overline{\Psi}%
{}^{(k,l)}(\eta)[H_{l,k+1} P_{l,k+1}(f)].
\]
This clearly yields that
\begin{eqnarray*}
\frac{\overline{\Phi}{}^{(k,l)}_{p}(\eta)(Q_{k+1}(f))}{\overline
{\Phi}%
{}^{(k,l)}_{p}(\eta)(Q_{k+1}(1))}&=&\frac{\overline{\Psi}{}^{(k,l)}_{p}(
\eta)[H_{l,k+1} P_{l,k+1}(f)]}{\overline{\Psi}{}^{(k,l)}_{p}(\eta)[H_{l,k+1}]}
\\
&=&\Psi{}^{H_{l,k+1}}_{l} \bigl( \overline{\Psi}{}^{(k,l)}_{p}(\eta
) \bigr)
P_{l,k+1}(f)
\end{eqnarray*}
and therefore
\begin{eqnarray*}
\overline{\Psi}{}^{(k+1),Q_{k+1}(1)}_{n}\bigl(\overline{\Phi}{}^{(k,l)}(\eta
)\bigr)P_{k+1,k+1} &=& \frac{1}{n+1}\sum_{p=0}^{n}\Psi
_{l}^{H_{l,k+1}} \bigl(
\overline{\Psi}{}^{(k,l)}_{p}(\eta) \bigr) P_{l,k+1}(f) \\
&=& \overline{\Psi}{}^{(l),H_{l,k+1}}_{n} \bigl( \overline{\Psi}%
{}^{(k,l)}(\eta) \bigr) P_{l,k+1}(f).
\end{eqnarray*}
In summary, we have proved that
\begin{eqnarray}
\overline{\Phi}{}^{(k+1,l)}(\eta)=\overline{\Psi}{}^{(k+1,l)}(\eta
)P_{l,k+1}(f)\hspace*{80pt}\nonumber\\
\eqntext{\mbox{with } \overline{\Psi}{}^{(k+1,l)}(\eta
)=\overline{\Psi}%
{}^{(l),H_{l,k+1}}_{n} \bigl( \overline{\Psi}{}^{(k,l)}(\eta) \bigr).}
\end{eqnarray}
This ends the proof of the proposition.
\end{pf}

\subsection{Contraction inequalities}

\begin{prop}
For any $l\leq k$ we have
\[
\beta(P_{l,k})=\frac{1}{2} \sup_{\eta,\mu}\bigl\Vert\overline{\Phi}{}^{(k,l)}(\eta)-%
\overline{\Phi}{}^{(k,l)}(\mu)\bigr\Vert.
\]
\end{prop}
\begin{pf}
Using Proposition \ref{laprop}, we find that
\begin{eqnarray*}
\bigl\Vert\overline{\Phi}{}^{(k,l)}(\eta)-\overline{\Phi}{}^{(k,l)}(\mu
)\bigr\Vert
&=&\bigl\Vert\bigl[ \overline{\Psi}{}^{(k,l)}(\eta)-\overline{\Psi
}{}^{(k,l)}(\mu
) \bigr] P_{l,k}\bigr\Vert\\
&\leq&\beta(P_{l,k}) \bigl\Vert\overline{\Psi}{}^{(k,l)}(\eta)-\overline
{\Psi}{}
^{(k,l)}(\mu)\bigr\Vert.
\end{eqnarray*}
This implies that
\[
\sup_{\eta,\mu}\bigl\Vert\overline{\Phi}{}^{(k,l)}(\eta)-\overline
{\Phi}%
{}^{(k,l)}(\mu)\bigr\Vert\leq2 \beta(P_{l,k}).
\]
On the other hand, if we chose the constant Dirac distribution flows
$\eta
=(\eta_{n})_{n\geq0}$ and $\mu=(\mu_{n})_{n\geq0}$ given by
\[
\forall n\geq0\qquad \eta_{n}=\delta_{x} \quad\mbox{and}\quad
\mu
_{n}=\delta_{y}
\]
for some $x,y\in S^{(l-1)}$, we also have that
\[
\overline{\Phi}{}^{(k,l)}(\delta_{x})-\overline{\Phi}{}^{(k,l)}(\delta
_{y})=\delta_{x}P_{l,k}-\delta_{y}P_{l,k}.
\]
This implies that
\[
\sup_{\eta,\mu}\bigl\Vert\overline{\Phi}{}^{(k,l)}(\eta)-\overline
{\Phi}%
{}^{(k,l)}(\mu)\bigr\Vert\geq{\sup_{x,y}}\Vert\delta_{x}P_{l,k}-\delta
_{y}P_{l,k}\Vert=2 \beta(P_{l,k}).
\]
This ends the proof of the proposition.
\end{pf}

Our next objective is to estimate the contraction coefficient $\beta
(P_{l,k})$ in terms of the mixing type properties of the semigroup $%
L_{l,k}=L_{l}L_{l-1}\cdots L_{k}$ associated with the Markov operators $
L_{l} $. We introduce the following regularity conditions.\vspace*{8pt}

$(L)_{m}$ \textit{There exists an integer $m\geq
1$ and a
sequence $(\varepsilon_{l}(L))_{l\geq0}\in(0,1)^{\mathbb{N}}$ such that}
\[
\forall l\geq0, \forall(x,y)\in\bigl(S^{(l-1)}\bigr)^{2}\qquad
L_{l+1,l+m}(x,%
\cdot)\geq\varepsilon_{l}(L) L_{l+1,l+m}(y,\cdot).
\]
\vspace*{0pt}

It is well known that the above condition is satisfied for any
aperiodic and
irreducible Markov chain on a finite space. Loosely speaking, for
noncompact spaces this condition is related to the tails of the transition
distributions on the boundaries of the state space. For instance, let us
assume that $S^{(l)}=\mathbb{R}$ and $L_{l}$ is the bi-Laplace transition
given by
\[
L_{l}(x,dy)=\frac{c(l)}{2} e^{-c(l) |y-A_{l}(x)|} \,dy
\]
for some $c(l)>0$ and some drift function $A_{n}$ with bounded
oscillations $\operatorname{osc}(A_{l})<\infty$. In this case, it is
readily checked that
condition $(L)_{m}$ holds true for $m=1$ with the parameter
\[
\varepsilon_{l-1}(L)=\exp{(-c(l) \operatorname{osc}(A_{l}))}.
\]
Under the condition (G) presented on page 11 and the mixing
condition $(L)_{m}$ stated above, we proved in \cite{fk} (see Corollary
4.3.3 on page 141) that we have for any $k\geq m\geq1$, and $l\geq1$
\[
\beta(P_{l+1,l+k})\leq\prod_{i=0}^{\lfloor k/m\rfloor-1} \bigl(
1-\varepsilon
_{l+im}^{(m)} \bigr)\qquad \mbox{with }
\varepsilon_{l}^{(m)}:=\varepsilon_{l}^{2}(L) \prod_{l+1\leq
k<l+m}\varepsilon_{k}(G).
\]
Several contraction inequalities can be deduced from these estimates, we
refer to Chapter 4 of the book \cite{fk}. To give a flavor of these results,
we further assume that $(M)_{m}$ is satisfied with $m=1$ and
$\varepsilon
(L)=\inf_{l}{\varepsilon_{l}(L)}>0$. In this case, we can check that
\[
\beta(P_{l+1,l+k})\leq\bigl( 1-\varepsilon(L)^{2} \bigr) ^{k}.
\]

\section*{Acknowledgments}

We would like to thanks the anonymous referee for helpful comments and
valuable suggestions that have improved the presentation of the article.

%

%
\printaddresses

\end{document}